\documentclass{article}%
\usepackage{latexsym}
\usepackage{amsmath}
\usepackage{graphicx}
\usepackage{amsfonts}
\usepackage{amssymb}%
\newtheorem{theorem}{Theorem}

\newtheorem{corollary}[theorem]{Corollary}

\newtheorem{definition}[theorem]{Definition}

\newtheorem{lemma}[theorem]{Lemma}
\newtheorem{notation}[theorem]{Notation}

\newtheorem{proposition}[theorem]{Proposition}
\newtheorem{remark}[theorem]{Remark}

\begin{document}

\author{Andrey Todorov\\UC Santa Cruz\\Department of Mathematics\\Santa Cruz, CA 95064\\Institute of Mathematics\\Bulgarian Academy of Sciences\\Sofia, Bulgaria}
\title{Mirror Symmetry, Borcherd-Harvey-Moore Products and \\Determinants of the Calabi-Yau Metrics on K3 Surfaces }
\maketitle
\date{}

\begin{abstract}
In the study of moduli of elliptic curves the Dedekind eta function
$\eta=q^{1/24}%
{\displaystyle\prod\limits_{n=1}^{\infty}}
\left(  1-q^{n}\right)  ,$ where $q=e^{2\pi i\tau}$ plays a very important
role. We will point out the three main properties of $\eta.$ \textbf{1. }It is
well known fact that $\eta^{24}$ is an automorphic form which vanishes at the
cusp. In fact $\eta^{24}$ is the discriminant of the elliptic curve.
\textbf{2. }The Kronecker limit formula gives the explicit relations between
the regularized determinant of the flat metric on the elliptic and $\eta.$
\textbf{3. }The Fourier expansion of $\frac{d}{dt}\log\eta(it)$ are positive
integers which give the number of elliptic curve that that are covering of the
elliptic \ curve $E_{\tau}$ of degree $n.$

Based on the work of Borcherds we construct on the moduli space of K3 surfaces
with B-field an automorphic form $\exp\Phi_{4,20}$ which vanishes on the
totally geodesic subspaces orthogonal to $-2$ vectors of $\mathbb{U}^{4}%
\oplus(\mathbb{E}_{8}(-1))^{2}$. We give an explicit formula of the
regularized determinants of the Laplacians of Calabi Yau metrics on K3
Surfaces, following suggestions by R. Borcherds. The holomorphic part of the
regularized determinants will be the higher dimensional analogue of Dedekind
Eta function.

We give explicit formulas for the number of non singular rational curves with
a fixed volume with respect to a Hodge metric in the case of K3 surfaces with
Picard group unimodular even lattice. The counting of rational curves on
special K3 surfaces using the regularized determinants of the Laplacian of CY
metrics is related to some results of Bershadsky, Cecotti, Ouguri and Vafa.
See \cite{BCOV}.

\end{abstract}
\tableofcontents

\section{Introduction}

\subsection{General Remarks}

In the study of moduli of elliptic curves the Dedekind eta function%
\[
\eta=q^{1/24}%
{\displaystyle\prod\limits_{n=1}^{\infty}}
\left(  1-q^{n}\right)  ,
\]
where $q=e^{2\pi i\tau}$ plays a very important role. We will point out the
three main properties of $\eta.$ \textbf{1.}It is well known fact that
$\eta^{24}$ is a automorphic form which vanishes at the cusp. In fact
$\eta^{24}$ is the discriminant of the elliptic curve. \textbf{2. }The
Kronecker limit formula gives the explicit relations between the regularized
determinant of the flat metric on the elliptic and $\eta.$ \textbf{3. }The
Fourier expansion of $\frac{d}{dt}\log\eta(it)$ are positive integers which
give the number of elliptic curve that that are covering of the elliptic
\ curve $E_{\tau}$ of degree $n.$

In this paper we will give the analogue of the Dedekind eta function for K3
surfaces. We will show that the main properties stated above of $\eta$ are
satisfied in the case of K3 surfaces.

The study of the moduli space of K3 surfaces recently attracted the attention
of string theorists. It is interesting that studies in optics by Frenel and
Hamilton was the first reason to study K3 surfaced. It was A. Weil who outline
the main problems in the study of the moduli of K3 surfaces. See \cite{W}. The
first main result in the study of moduli of K3 surfaces is due to Shafarevich
and Piatetski-Shapiro. See \cite{PS}. They proved the global Torelli Theorem
for polarized algebraic K3 surfaces. Combining the Theorem of Shafarevich and
Piatetski Shapiro with the description of the mapping class group of K3
surface one obtain that the moduli space $\mathfrak{M}_{k3,n}$ of polarized
algebraic K3 surfaced with a polarization class $e$ such that $\left\langle
e,e\right\rangle =2n>0$ is a Zariski open set in
\[
\Gamma_{K3,n}^{+}\backslash\mathbb{SO}(2,19)/\mathbb{SO}(2)\times
\mathbb{SO}(19),
\]
where $\Gamma_{K3,n}^{+}$ is an index two subgroup in $\mathcal{O}%
_{\Lambda_{K3,n}}(\mathbb{Z})$ and $\Lambda_{K3,n}$ is the lattice isomorphic
to $-2n\mathbb{Z\oplus U}^{2}\mathbf{\oplus}\mathbb{E}_{8}(-1)\oplus
\mathbb{E}_{8}(-1).$ In \cite{To80} it was proved that every point of
$\mathbb{SO}(3,19)/\mathbb{SO}(2)\times\mathbb{SO}(1,19)$ corresponds to a
marked K3 surface. Based on this result in \cite{KT} it was proved that the
moduli space of Ricci flat metrics on K3 surfaces with a fixed volume is
isomorphic to
\[
\mathfrak{M}_{KE}:=\Gamma^{+}\backslash\left(  \mathbb{SO}_{0}%
(3,19)/\mathbb{SO}(3)\times\mathbb{SO}(19)-\mathcal{D}_{KE}\right)  ,
\]
where $\Gamma^{+}$ is a subgroup of index 2 in the group of automorphisms of
the group of the automorphisms of the Euclidean lattice $\Lambda
_{K3}=\mathbb{U}^{3}\mathbf{\oplus}\mathbb{E}_{8}(-1)\oplus\mathbb{E}%
_{8}(-1),$ where%
\[
\mathbb{U=}\left(
\begin{array}
[c]{cc}%
0 & 1\\
1 & 0
\end{array}
\right)
\]
is the hyperbolic lattice and $\mathbb{E}_{8}(-1)$ is the standard lattice and
$\mathcal{D}_{KE}$ is the subspace whose points correspond to Ricci flat
metrics on orbifolds. Donaldson proved in \cite{D} that the mapping class
group $\Gamma$ of a K3 surface is\ a subgroup of index 2 in the group of the
automorphisms of the Euclidean lattice $\Lambda_{K3}$.

Mirror Symmetry is based on the observation that there are two different
models A and B in string theory whcih define one and the same partition
function. The A model is related to the deformation of the K\"{a}hler-Einstein
metrics. The B-model is related to the deformations of complex structures. To
studedy mirror symmetryc on K3 surfaces we need to define a B-field on a K3
surfaces. It is a class of cohomology $\omega_{X}(1,1)\in H^{1,1}%
($X$,\mathbb{C})$ of type $(1,1)$ on a K3 surface X such that
\[%
{\displaystyle\int\limits_{\text{X}}}
\operatorname{Im}\omega\wedge\operatorname{Im}\omega>0.
\]
The moduli space of marked K3 surfaces with a B-field is isomorphic to
$\mathfrak{h}_{4,20}:=\mathbb{SO}_{0}(4,20)/\mathbb{SO}(4)\times
\mathbb{SO}(20).$ Aspinwall and Morrison proved that the moduli space of Super
Conformal Field Theories with supersymmetry (4,4) is described by $\Gamma
_{B}^{+}\backslash\mathfrak{h}_{4,20},$ where $\Gamma_{B}^{+}$ is a subgroup
of index two in $\mathcal{O}(\Lambda_{K3}).$ It is well known that
$\mathfrak{h}_{4,20}$ parametrizes the four-dimensional oriented subspaces in
$\mathbb{R}^{4,20}$ on which the bilinear form is strictly positive. See
\cite{AM}. To a pair (X,$\omega_{X}(1,1))$ of a marked K3 surface with a
B-field $\omega_{X}(1,1)$ we assign a oriented four dimensional subspace
$E_{\text{X,}\omega_{X}(1,1)}$ in
\[
H^{\ast}(\text{X,}\mathbb{Z})\otimes\mathbb{R}=\left(  H^{0}(\text{X,}%
\mathbb{Z})\oplus H^{2}(\text{X,}\mathbb{Z})\oplus H^{4}(\text{X,}%
\mathbb{Z})\right)  \otimes\mathbb{R}%
\]
on which the bilinear form defined by the cup product is positive. We will
assume that $\left(  H^{0}(\text{X,}\mathbb{Z})\oplus H^{4}(\text{X,}%
\mathbb{Z})\right)  =\mathbb{U}_{0}$ and the B-field $\omega_{X}(1,1)$ we will
be identified with
\begin{equation}
(1,-\frac{1}{2}\left(  \omega_{X}(1,1)\wedge\omega_{X}(1,1)\right)
)\in\left(  H^{0}(\text{X,}\mathbb{Z})\oplus H^{2}(\text{X,}\mathbb{Z})\oplus
H^{4}(\text{X,}\mathbb{Z})\right)  . \label{ext}%
\end{equation}
From now on we will consider the B-field $\omega_{X}(1,1)$ as defined by
$\left(  \ref{ext}\right)  .$ The four dimensional subspace $E_{\text{X,}%
\omega_{X}(1,1)}$ contains the two dimensional subspace $E_{\omega_{\text{X}}%
}$ spanned by $\operatorname{Re}\omega_{\text{X}}$ and $\operatorname{Im}%
\omega_{\text{X}},$ where $\omega_{\text{X}}$ is the holomorphic two form on X
defined up to a constant and the two dimensional subspace $E_{\omega_{X}%
(1,1)}$ spanned by $\operatorname{Re}\omega_{X}(1,1)$ and $\operatorname{Im}%
\omega_{X}(1,1),$ where $\omega_{X}(1,1)$ is defined by $\left(
\ref{ext}\right)  $. $E_{\omega_{\text{X}}}$ will the orthogonal to
$E_{\omega_{\text{X}}(1,1)}$ in $E_{\text{X,}\omega_{X}(1,1)}.$

Mirror Symmetry is pretty well understood in the case of K3 surfaces. See
\cite{AM}, \cite{dol} and \cite{To93}. The mirror symmetry is exchanging
$E_{\omega_{\text{X}}}$ with $E_{\omega_{\text{X}}(1,1)}.$ Special case of
mirror symmetry of algebraic K3 surfaces was studied in details in \cite{dol}.

In this paper we will consider the moduli space of K3 surfaces with $B$-fields
and an automorphic form $\exp\left(  \Phi_{4,20}\right)  $ which vanishes on
the totally geodesic subspaces that are orthogonal to $-2$ vectors. Such
automorphic form exists according to \cite{B97}. The restriction of
$\exp\left(  \Phi_{4,20}\right)  $ on the moduli space \ $\mathfrak{M}%
_{ell}:=\Gamma_{ell}\backslash\mathfrak{h}_{2,10}$ of elliptic K3 surfaces
with a section vanishes on the discriminant locus $\mathfrak{D}_{ell}%
\subset\mathfrak{M}_{ell}\subset\Gamma_{B}^{+}\backslash\mathfrak{h}_{4,20}$
which is defined by the points orthogonal to $-2$ vectors. The mirror $Y$ of
the elliptic K3 $X$ with the section has Picard group$Pic(Y)=\mathbb{U}%
\oplus\mathbb{E}(-1)\oplus\mathbb{E}(-1).$ $\exp\left(  \Phi_{4,20}\right)  $
restricted on a line $tL$ in the K\"{a}hler cone $K(Y)$ spanned by the
imaginary part $L$ of a Hodge metric, has a Fourier expansion. The Fourier
coefficients of $\frac{d}{dt}\log\Phi(it)$ are positive integers and they
count the number of rational curves of fix volume. Thus in the A model of K3
surfaces with Picard group unimodular lattice, the restriction of $\exp\left(
\Phi_{4,20}\right)  $ on the K\"{a}hler cone counts rational curves. In the
B-model the restriction of $\exp\left(  \Phi_{4,20}\right)  $ on the moduli
space counts vanishing invariant cycles.

The regularized determinants of the Laplacian of Ricci flat metrics
$\det(\Delta_{KE})$ acting on $(0,1)$ forms will be a function on on the
moduli space of Einstein metric $\mathfrak{M}_{KE}.$ R. Borcherds suggested
that one can compute the determinants of the Laplacians of Ricci flat metrics
explicitly by using the method of the theta lifts. See \cite{B97}. In this
paper we will give an explicit expression of the regularized determinants of
the Laplacians of CY metrics $\det$ as a function on the moduli space of
Einstein metrics $\mathfrak{M}_{KE}.$

There are some relations of this paper with the papers \cite{CD} and \cite{CM}.

\subsection{Organization of the Paper}

In \textbf{Section 2 }we describe some basic property of the symmetric space
\[
\mathfrak{h}_{p,q}:=\mathbb{SO}_{0}(p,q)/\mathbb{SO}(p)\times\mathbb{SO}(q).
\]

In \textbf{Section 3 }we study the unimodular even indefinite lattices
$\Lambda_{p,q}.$ We define the discriminant locus $\mathfrak{D}_{p,q}$ in the
locally symmetric space $\mathcal{O}_{+}^{\ast}(\Lambda_{p,q})\backslash
\mathfrak{h}_{p,q}.$ We prove in this \textbf{Section }that $\mathfrak{D}%
_{p,q}$ is irreducible.

In \textbf{Section 4} we describe the main results about moduli of K3 surfaces.

In \textbf{Section 5 }we study automorphic forms on $\mathcal{O}_{+}^{\ast
}(\Lambda_{p,q})\backslash\mathfrak{h}_{p,q}.$

In \textbf{Section 6} we prove the analogue of the Kronecker limit formula,
i.e. we gave the explicit formula for the determinant of the Laplacian of
Calabi-Yau metrics (K\"{a}hler-Einstein metrics) on the moduli space of
K\"{a}hler-Einstein metrics $\mathcal{O}_{+}^{\ast}(\Lambda_{3,19}%
)\backslash\mathfrak{h}_{3,19}.$

In \textbf{Section 7 }we study mirror symmetry of K3 surfaces.

In \textbf{Section 8 }we constructed the analogue of the Dedekind eta function
for K3 surfaces and proved its main properties.

\subsection{Acknowledgements}

The author want to acknowledge the help and suggestions of Greg Zuckerman. He
proposed to study the behavior of the regularized determinants on the moduli
space of elliptic curves with sections.

Special thanks to Jay Jorgenson for his help and comments. I am grateful to
him for introducing me to the ideas of regularized determinants.

I want to thank Jun Li for his interest in this paper and help. Special thanks
to the Center of Mathematical Sciences of Zhe Jiang University and National
Center for Theoretical Sciences, Mathematical Division, National Tsing Hua
University for the financial support during the preparation of the paper.

\section{Symmetric Space $\mathfrak{h}_{p,q}:=\mathbb{SO}_{0}(p,q)/\mathbb{SO}%
(p)\times\mathbb{SO}(q)$}

\subsection{Global Flat Coordinates on the Symmetric Space $\mathfrak{h}%
_{p,q}$}

We will need some basic facts about the symmetric space
\[
\mathfrak{h}_{p,q}:=\mathbb{SO}_{0}(p,q)/\mathbb{SO}(p)\times\mathbb{SO}(q).
\]
The following Theorem is standard.

\begin{theorem}
\label{G}There is a one to one correspondence between points $\tau$ in
$\mathfrak{h}_{3,19}$ and all oriented three dimensional $E_{\tau}$ subspaces
in $\Lambda_{K3}\otimes\mathbb{R}$ on which the intersection form on
$\Lambda_{K3}$ is strictly positive.
\end{theorem}

\begin{theorem}
\label{G4}Let $\mathbb{R}^{p,q}$ be a $p+q$ dimensional real vector space with
a metric with signature $(p,q).$ Let $E_{\tau_{0}}$ be a $p-$dimensional
subspace in $\mathbb{R}^{p,q}$ such the restriction of the quadratic form on
$E_{\tau_{0}}$ is strictly positive. Let $e_{1},...e_{p}$ be an orthonormal
basis of $E_{\tau_{0}}.$ Let $e_{p+1},...,e_{p+q}$ be orthogonal vectors to
$E_{\tau_{0}}$ such that $\left\langle e_{i},e_{j}\right\rangle =-\delta_{ij}$
for $p+1\leq i,j\leq p+q.$ Let $E_{\tau}$ be any $p-$dimensional subspace in
$\mathbb{R}^{p,q}$ such that the restriction of the quadratic form in
$E_{\tau}$ is strictly positive$.$ Then there exists a basis $\{g_{1}%
(\tau),...,g_{p}(\tau)\}$ in $E_{\tau}$ such that
\begin{equation}
g_{j}(\tau)=e_{j}+\sum_{i=p+1}^{p+1}\tau_{j}^{i}e_{i}. \label{g1}%
\end{equation}

\end{theorem}

\textbf{Proof:} Let $f_{1}(\tau),$...,$f_{p}(\tau)$ be an orthonormal basis of
$E_{\tau}.$ Let
\begin{equation}
f_{i}=\sum_{j=p+1}^{p+q}\tau_{i}^{j}e_{j} \label{g4}%
\end{equation}
for $1\leq i\leq p$ and $1\leq j\leq p+q.$ Let $\left(  A_{ij}(\tau)\right)  $
be the $p\times p$ matrix $\left(  \tau_{i}^{j}\right)  $ for $1\leq i\leq p$
and $1\leq j\leq q$ where $\tau_{i}^{j}$ are the elements in the expression
$\left(  \ref{g4}\right)  .$

\begin{lemma}
\label{G41}$\det(A_{ij}(\tau))\neq0.$
\end{lemma}

\textbf{Proof:} Suppose that $\det(A_{ij}(\tau))=0.$ This implies that
$rk(A_{ij}(\tau))<p.$ So we can find constants $a_{i}$ for $i=1,...,p$ such
that at least one of them is non zero and
\begin{equation}%
{\displaystyle\sum\limits_{j=1}^{p}}
a_{j}\left(
{\displaystyle\sum\limits_{i=1}^{p}}
\tau_{i}^{j}\right)  =0. \label{g5}%
\end{equation}
Let us consider the vector
\begin{equation}
g(\tau)=%
{\displaystyle\sum\limits_{j=1}^{p}}
a_{j}g_{j}. \label{g6}%
\end{equation}
Combining $\left(  \ref{g5}\right)  $ and $\left(  \ref{g6}\right)  $ we
obtain that
\begin{equation}
g(\tau)=\sum_{i=p+1}^{p+q}\lambda_{i}e_{i}. \label{g7}%
\end{equation}
The definition of the vectors $g_{i}(\tau)$ and $\left(  \ref{g6}\right)  $
imply that
\begin{equation}
\left\langle g(\tau),g(\tau)\right\rangle =-2\sum_{i=p+1}^{p+q}\left\vert
\lambda_{i}\right\vert ^{2}<0. \label{g9}%
\end{equation}
Clearly $g(\tau)$ is a non zero vector in $E_{\tau}.$ So $\left\langle
g(\tau),g(\tau)\right\rangle >0.$ Thus we get a contradiction with $\left(
\ref{g9}\right)  .$ Lemma \ref{G41} is proved. $\blacksquare.$

Theorem \ref{G4} follows directly from Lemma \ref{G41}. $\blacksquare.$

\begin{corollary}
\label{G43}\textit{There is one to one correspondence between the set of all
}$p\times q$\textit{\ matrices (}$\tau_{i}^{j})$ for $1\leq i\leq p$ and
$p+1\leq j\leq p+q$\textit{\ such that the vectors }$g_{i}(\tau)$\textit{\ for
}$i=1,...,p$ \textit{defined by }$\left(  \ref{g1}\right)  $\textit{\ spanned
a }$p-$\textit{dimensional subspace }$E_{\tau}$\textit{\ in }$\mathbb{R}%
^{p,q}$\textit{\ on which the restriction of the quadratic form }$\left\langle
u,v\right\rangle $\textit{\ is strictly positive and the set of points in
}$\mathfrak{h}_{3,19}$\textit{. Thus }$(\tau_{i}^{j})$ define global
coordinates on $\mathfrak{h}_{3,19}.$
\end{corollary}

\subsection{Decomposition of $\mathfrak{h}_{p,q}$}

\begin{theorem}
\label{Dec1}We have the following decomposition of
\begin{equation}
\mathfrak{h}_{2,p}=\mathbb{R}^{1,p-1}\mathfrak{h}_{1,p-1}+\sqrt{-1}%
\mathfrak{h}_{1,p-1}. \label{F1}%
\end{equation}

\end{theorem}

\textbf{Proof: }It is a well known fact that $\mathfrak{h}_{1,p-1}$ is one of
the component $V^{+}$ of the cone $V:=\left\{  v\in\mathbb{R}^{1,p-1}%
|\left\langle v,v\right\rangle >0\right\}  .$ Let us consider $\mathbb{R}%
^{2,p}=\mathbb{R}^{1,p-1}\oplus\mathbb{R}^{1,1}.$ Let us consider the map:
\[
w\in\mathbb{R}^{1,p-1}\mathfrak{h}_{1,p-1}+\sqrt{-1}\mathfrak{h}%
_{1,p-1}\rightarrow\mathbb{P}\left(  \left(  \mathbb{R}^{1,p-1}\oplus
\mathbb{R}^{1,1}\right)  \otimes\mathbb{C}\right)
\]
defined as follows%
\[
\Psi:w=(w_{1},...,w_{p})\rightarrow\left(  w_{1},...,w_{p},-\frac{\left\langle
w,w\right\rangle }{2},1\right)  .
\]
It is easy to check that in $\mathbb{P}\left(  \mathbb{R}^{2,p}\otimes
\mathbb{C}\right)  $ we have $\left\langle \Psi(w),\Psi(w)\right\rangle =0$
$and$ $\left\langle \Psi(w),\overline{\Psi(w)}\right\rangle >0.$ Thus the
image of $\mathbb{R}^{1,p-1}\mathfrak{h}_{1,p-1}+\sqrt{-1}\mathfrak{h}%
_{1,p-1}$ under the map $\Psi$ will be $\mathfrak{h}_{2,p},$ since
$\mathfrak{h}_{2,p}$ in $\mathbb{P}\left(  \mathbb{R}^{2,p}\otimes
\mathbb{C}\right)  $ is given by one of the components of the open set in the
quadratic $\left\langle w,w\right\rangle =0$ defined by $\left\langle
w,\overline{w}\right\rangle >0.$ It is very easy to prove that $\Psi$ is one
to one map. $\blacksquare$

\begin{theorem}
\label{Dec}Suppose that $p\geq3,$ and $q\geq2.$ Then we have the following
decomposition
\begin{equation}
\mathfrak{h}_{p,q}=\mathfrak{h}_{p-1,q-1}\times\mathbb{R}^{p-1,q-1}%
\times\mathbb{R}_{+}, \label{F}%
\end{equation}
where $\mathbb{R}_{+}$ is the set of real positive numbers$.$
\end{theorem}

\textbf{Proof: }Let us consider in the space $\mathbb{R}^{p,q}$ two vectors
$e_{p+q-1}$ and $e_{p+q}$ such that%
\[
\left\langle e_{p+q},e_{p+q}\right\rangle =\left\langle e_{p+q-1}%
,e_{p+q-1}\right\rangle =0\text{ and }\left\langle e_{p+q-1},e_{p+q}%
\right\rangle =1.
\]
Clearly the orthogonal complement to the subspace $\{e_{p+q},e_{p+q}\}$ will
be isometric to $\mathbb{R}^{p-1,q-1}.$ Let us consider a basis $\{e_{1}%
,...,e_{p+1}\}$ of $\mathbb{R}^{p,q},$ where $e_{1},...,e_{p+q-2}$ is a basis
of $\mathbb{R}^{p-1,q-1}.$

There is one to one correspondence between the points $\tau\in\mathfrak{h}%
_{p,q}$ and the oriented $p-$dimensional subspaces $E_{\tau}$ in
$\mathbb{R}^{p,q}$ on which the restriction of the bilinear form is strictly
positive. The intersection $E_{\tau}\cap\mathbb{R}^{p-1,q-1}$ will be $\left(
p-1\right)  -$dimensional subspace in $\mathbb{R}^{p-1,q-1}$ on which the
bilinear form is strictly positive. Let $f_{1}$ be a vector in $\mathbb{R}%
^{p,q}$ orthogonal to $\mathbb{R}^{p-1,q-1}\cap E_{\tau}.$ It is easy to see
that the coordinates of $f_{1}$ can be normalized in such a way that its
coordinates in $\mathbb{R}^{p,q}$ are $f_{1}=(\mu_{1},...,\mu_{20}%
,1,\lambda),$ where $\mu=(\mu_{1},...,\mu_{p+q-2})$ is any vector in
$\mathbb{R}^{p-1,q-1}$ and $\lambda>0$ and $\lambda>\left\langle \mu
,\mu\right\rangle .$ Thus the correspondence $E_{\tau}\rightarrow\left(
f_{1},E_{\tau}\cap\mathbb{R}^{p-1,q-1}\right)  $ establishes the decomposition
$\left(  \ref{F}\right)  .$ $\blacksquare$

\subsection{Definition of the Standard Metric on $\mathfrak{h}_{p,q}$}

Since $\mathfrak{h}_{p,q}\subset Grass(p,p+q)$ then the tangent space
$T_{\tau_{0},\mathfrak{h}_{p,q}}$ at a point $\tau_{0}\in\mathfrak{h}_{p,q}$
can be identified with $Hom\left(  E_{\tau_{0}},E_{\tau_{0}}^{\perp}\right)
.$ Thus any tangent vector $A\in T_{\tau_{0},\mathfrak{h}_{p,q}}$ can be
written in the form%
\begin{equation}
A=%
{\displaystyle\sum\limits_{i=1}^{p}}
{\displaystyle\sum\limits_{j=p+1}^{p+q}}
\tau_{i}^{j}\left(  e_{i}^{\ast}\otimes e_{j}\right)  , \label{M}%
\end{equation}
where $e_{i}$ for $i=1,...,q$ is an orthonormal basis of $E_{\tau_{0}}$ and
$e_{j}$ for $j=p+1,...,p+q$ is an orthonormal basis of \ $E_{\tau_{0}}^{\perp
}.$ Then we define the metric on $T_{\tau_{0},\mathfrak{h}_{p,q}}=Hom\left(
E_{\tau_{0}},E_{\tau_{0}}^{\perp}\right)  $ for $A\in T_{\tau_{0}%
,\mathfrak{h}_{p,q}}=Hom\left(  E_{\tau_{0}},E_{\tau_{0}}^{\perp}\right)  $
given by $\left(  \ref{M}\right)  $ by
\begin{equation}
\left\Vert A^{2}\right\Vert =%
{\displaystyle\sum\limits_{i,j}}
\left\vert \tau_{i}^{j}\right\vert ^{2}. \label{met}%
\end{equation}
We will call this metric the Bergman metric on $\mathfrak{h}_{3,19}.$

\begin{lemma}
\label{MET}The Bergman metric $ds_{B}^{2}$ is invariant metric on
$\mathfrak{h}_{p,q}.$ It is given in the flat coordinate system by%
\begin{equation}
ds_{B}^{2}=%
{\displaystyle\sum\limits_{1\leq j\leq3,1\leq i19}}
\left(  d\tau_{j}^{i}\right)  ^{2}+O(2). \label{met1}%
\end{equation}

\end{lemma}

\textbf{Proof: }The proof of Lemma \ref{MET} follows directly from the
definition of the Bergman metric. $\blacksquare$

\section{Discriminants in $\mathfrak{H}_{p,q}$}

\subsection{Definition and Basic Properties of the Discriminant}

From now on we will consider the symmetric spaces $\mathfrak{h}_{p,q}$ for
which $p-q\equiv0\operatorname{mod}8.$ From now on this paper $\Lambda p,q$
will be unimodular even lattice of signature $(q+8k,q).$ We have the following
description all $\Lambda p,q:$

\begin{theorem}
\label{uel}Suppose that $\Lambda p,q$ the unimodular even lattice of signature
$(p,q)$ for $p-q\equiv0\operatorname{mod}8.$ Then
\[
\Lambda_{p,q}\approxeq\underset{p=q+8k}{\underbrace{\mathbb{U}\oplus
...\oplus\mathbb{U}}}\oplus\underset{q}{\underbrace{\mathbb{E}_{8}%
(-1)\oplus...\oplus\mathbb{E}_{8}(-1)}}.
\]

\end{theorem}

\begin{definition}
Define the set $\Delta_{p,q}(e):=\{\delta\in\Lambda_{p,q}|$ $\left\langle
\delta,\delta\right\rangle =-2\}.$ Let us define by $\mathcal{O}_{p,q}$ the
group of the automorphisms of the lattice $\Lambda_{p,q}.$ Let $\mathcal{O}%
_{p,q}^{+}$ be the subgroup of $\mathcal{O}_{p,q}$ which preserve the
orientation of the positive subspaces of dimension $p$ in $\Lambda
_{p,q}\otimes\mathbb{R}.$ Then $\mathcal{O}_{p,q}^{+}$ has index two in
$\mathcal{O}_{p,q}.$
\end{definition}

\begin{definition}
We know that $\mathfrak{h}_{p,q}$ can be realized as an open set in the
Grassmanian $Grass(p,p+q).$ Let us denote by $\mathfrak{h}_{p,q-1}(\delta)$
the set of all $p-$dimensional subspaces in the orthogonal complement of the
vector $\delta$ in $\Lambda_{p,q}\otimes\mathbb{R}.$ We will define the
discriminant locus $\mathfrak{D}_{p,q}$ in $\mathcal{O}_{p,q}^{+}%
\backslash\mathfrak{h}_{p,q}$ as follows:%
\[
\mathfrak{D}_{p,q}:=\mathcal{O}_{p,q}^{+}\backslash\left(  \underset{\delta
\in\Delta(e)}{\cup}(\mathfrak{h}_{p,q-1}(\delta))\right)
\]

\end{definition}

This definition is motivated by the definition of the discriminant locus in
the moduli of algebraic K3 surfaces.

\subsection{The Irreducibility of the Discriminant}

\begin{theorem}
\label{Bor}The discriminant locus $\mathcal{D}_{p,q}$ in $\mathcal{O}%
_{p,q}^{+}\backslash\mathfrak{h}_{p,q}$ is an irreducible divisor, where
$\Lambda_{p,q}$ is an even unimodular lattice.
\end{theorem}

\textbf{Proof:} The proof of Theorem \ref{Bor} will follow if we prove that on
the set of vectors $\Delta_{\Lambda_{p,q}}$ the group $\mathcal{O}%
_{\Lambda_{p,q}}^{+}$ acts transitively. Thus they form one orbit and
therefore the discriminant locus $\mathcal{D}_{p,q}$ in $\mathcal{O}_{p,q}%
^{+}\backslash\mathfrak{h}_{p,q}$ is an irreducible divisor.

The proof that on the set of vectors $\Delta_{\Lambda_{p,q}}$ the group
$\mathcal{O}_{\Lambda_{p,q}}^{+}$ acts transitively will use similar ideas as
the proof of the irreducibility of the discriminant locus in the moduli space
of Enriques surfaces given by Borcherds in \cite{B96}.

We will proceed by induction on $p$ to prove that the action of $\mathcal{O}%
_{\Lambda_{p,q}}^{+}$ on the set $\Delta_{\Lambda_{p,q}}$ is transitive. For
$p=0$ the Theorem \ref{Bor} is obvious. Suppose that Theorem \ref{Bor} is true
for $p.$ We will denote by $L$ the lattice
\[
\underset{p}{\underbrace{\mathbb{U}\oplus...\oplus\mathbb{U}}}\oplus
\underset{q}{\underbrace{\mathbb{E}_{8}(-1)\oplus...\oplus\mathbb{E}_{8}(-1)}}%
\]
and by $M$ the lattice $L\oplus\mathbb{U}$.

The plan of the proof is the following. We will denote by $R_{0}$ and $R_{1}$
the set of norm $-2$ vectors of $M$ which have inter product respectively $0$
or $1$ with the vector $e=(\overrightarrow{0},0,1)\in L\oplus\mathbb{U=}M.$
Let $\Gamma_{1}$ be the group generated by reflections of elements of the set
$R_{1}$ and $\Gamma_{2}$ be the group generated by reflections of elements of
$R_{0}\cup R_{1}$ and $-id.$ We will show first that any $-2$ vector of $M$ is
conjugate to an element of the set $R_{0}\cup R_{1}.$ Then we will show that
the group $\mathcal{O}_{M}^{+}(\mathbb{Z})$ interchange the sets $R_{0}$ and
$R_{1}.$

\begin{lemma}
\label{Bor3}Any norm $-2$ vector $\delta$ of $M$ is conjugate to an element of
$R_{0}\cup R_{1}$ under the group $\Gamma_{1}.$
\end{lemma}

\textbf{Proof:} The proof of Lemma \ref{Bor3} is based on the following
Propositions \ref{Bor1} and \ref{Bor2}:

\begin{proposition}
\label{Bor1}Suppose that $v\notin L\subset L\otimes\mathbb{Q}.$ Suppose that
$x$ is some real number. Then there exists a vector $\overrightarrow{\mu}\in
L$ such that$\left\vert \left(  \overrightarrow{v}-\overrightarrow{\mu
}\right)  ^{2}-x\right\vert <1.$
\end{proposition}

\textbf{Proof:} The proof of Proposition \ref{Bor1} follows the proof of Lemma
\textbf{2.1} in \cite{B96}. Since $\overrightarrow{v}\notin L\subset
L\otimes\mathbb{Q}$ we can find a primitive isotropic vector $\overrightarrow
{\rho}$ such that $\left\langle \overrightarrow{\rho},\overrightarrow
{v}\right\rangle $ is not an integer. This is because primitive isotropic
vectors span $L.$ As the group $\mathcal{O}_{L}(\mathbb{Z})$ acts transitively
on norm $0$ vectors we can assume that
\[
\rho=(0,0,1)\in\underset{p-1}{\underbrace{\mathbb{U}\oplus...\oplus\mathbb{U}%
}}\oplus\underset{q}{\underbrace{\mathbb{E}_{8}(-1)\oplus...\oplus
\mathbb{E}_{8}(-1)}}\oplus\mathbb{U}%
\]
Then $\overrightarrow{v}=(\overrightarrow{\lambda},a,b)$ with $a$ not an
integer. We will find some $\overrightarrow{\mu}$ of the form $\overrightarrow
{\mu}$ $=(\overrightarrow{0},m,n)$ with integers $m$ and $n$ such that
\[
\left\vert \left\langle \overrightarrow{\mu}-\overrightarrow{v}%
,\overrightarrow{\mu}-\overrightarrow{v}\right\rangle -x\right\vert
=\left\vert \left\langle \overrightarrow{v},\overrightarrow{v}\right\rangle
-2(a-m)(b-n)-x\right\vert <1.
\]
Since $a$ is not an integer we can find some integer $m$ such that $\left\vert
a-m\right\vert <1.$ Whenever we add $1$ to $n,$ the expression $2(a-m)(b-n)$
is changed by a non zero number less than $2,$ so we can choose some integer
$n$ such that $2(a-m)(b-n)$ is at a distance of less then $1$ from any given
number $x-\lambda^{2}.$ This proves Proposition \ref{Bor1}. $\blacksquare$

\begin{proposition}
\label{Bor2}Suppose that $R_{1}$ is the set of norm $-2$ vectors of $M$ having
inner product $1$ with $e=(\overrightarrow{0},0,1)\in L\oplus\mathbb{U}.$
Suppose that $\Gamma_{1}$ is the subgroup of $\mathcal{O}_{M}(\mathbb{Z})$
generated by reflections of vectors of $R_{1}$ and the automorphism $-id.$
Then any vector $r\in M$ is conjugate under $\Gamma_{1}$ to a vector of the
form $(\overrightarrow{v},m,n)\in M$ such that either $m=0$ or $\frac
{\overrightarrow{v}}{m}\in L$ and $m>0.$
\end{proposition}

\textbf{Proof:} We can assume that $\overrightarrow{r}=(\overrightarrow
{v},m,n)$ has the property that $\left\vert \left\langle \overrightarrow
{r},\overrightarrow{e}\right\rangle \right\vert =|m|$ is minimal among all
conjugates of $r$ under $\Gamma_{1},$ where $\overrightarrow{e}%
=(\overrightarrow{0},0,1)$. If $m=0$ then we are done. So we can assume that
$m\neq0,$ and wish to prove that $\frac{v}{m}\in L$ and $m>0.$ Suppose that
$\frac{\overrightarrow{v}}{m}\notin L.$ By Proposition \ref{Bor1} we can find
a vector $\overrightarrow{\mu}\in L$ satisfying%
\begin{equation}
\left\vert \left\langle \overrightarrow{\mu}-\frac{\overrightarrow{v}}%
{m},\overrightarrow{\mu}-\frac{\overrightarrow{v}}{m}\right\rangle +\left(
-\frac{2n}{m}-\frac{\left\langle \overrightarrow{v},\overrightarrow
{v}\right\rangle }{m}\right)  \right\vert <1. \label{b0}%
\end{equation}
Let $\delta=\left(  \overrightarrow{\mu},1,\frac{-\left\langle \overrightarrow
{\mu},\overrightarrow{\mu}\right\rangle -2}{2}\right)  .$ It is easy to see
that $\left\langle \delta,\delta\right\rangle =-2.$ Let $T_{\delta
}(r)=r^{\prime}=r+\left\langle r,\delta\right\rangle \delta,$ i.e. $r^{\prime
}$ is the reflection of $r$ with respect to the hyperplane of $\delta\in
R_{1}.$ Direct computations show that%
\[
\left\vert \left\langle r^{\prime},\overrightarrow{e}\right\rangle \right\vert
=\left\vert \left\langle T_{\delta}(r),\overrightarrow{e}\right\rangle
\right\vert =\left\vert \left\langle r,T_{\delta}(\overrightarrow
{e})\right\rangle \right\vert =\left\vert \left\langle r,\overrightarrow
{e}+\delta\right\rangle \right\vert =
\]%
\begin{equation}
\left\vert m\left(  \left\langle \overrightarrow{\mu}-\frac{\overrightarrow
{v}}{m},\overrightarrow{\mu}-\frac{\overrightarrow{v}}{m}\right\rangle
+\left(  -\frac{2n}{m}-\frac{\left\langle \overrightarrow{v},\overrightarrow
{v}\right\rangle }{m^{2}}\right)  \right)  \right\vert . \label{b2}%
\end{equation}
Combining $\left(  \ref{b0}\right)  $ and $\left(  \ref{b2}\right)  $ we
deduce that%
\begin{equation}
\left\vert \left\langle r^{\prime},\overrightarrow{e}\right\rangle \right\vert
=\left\vert m\left\langle \overrightarrow{\mu}-\frac{\overrightarrow{v}}%
{m},\overrightarrow{\mu}-\frac{\overrightarrow{v}}{m}\right\rangle +\left(
-\frac{2n}{m}-\frac{\left\langle \overrightarrow{v},\overrightarrow
{v}\right\rangle }{m}\right)  \right\vert <m. \label{b3}%
\end{equation}
We have chosen
\begin{subequations}
\begin{equation}
\left\vert \left\langle r,\delta\right\rangle \right\vert =m \label{b4}%
\end{equation}
to be minimal. So we $\left(  \ref{b3}\right)  $ contradicts $\left(
\ref{b4}\right)  $. Proposition \ref{Bor2} is proved. $\blacksquare$

\textbf{Proof of Lemma }\ref{Bor3}: Let $\delta=(v,m,n).$ By Proposition
\ref{Bor2} we can assume that either $m=0$ or $\frac{\overrightarrow{v}}{m}\in
L.$ If $m=0$ then Lemma \ref{Bor3} is proved. Suppose that $\frac
{\overrightarrow{v}}{m}\in L$ holds. Then $\left\langle \frac{\overrightarrow
{v}}{m},\frac{\overrightarrow{v}}{m}\right\rangle \in\mathbb{Z}.$ Thus
\end{subequations}
\begin{equation}
\left\langle \delta,\delta\right\rangle =-2=m^{2}\left\langle \frac
{\overrightarrow{v}}{m},\frac{\overrightarrow{v}}{m}\right\rangle +2mn.
\label{b7}%
\end{equation}
So $\left(  \ref{b7}\right)  $ implies that $-2$ is divisible by $(m^{2},2m).$
From here we conclude that $m=1.$ Lemma \ref{b3} is proved. $\blacksquare$

Let us define the group $\Gamma_{3}$ as the group generated by the
automorphisms $\mathcal{O}_{L}(\mathbb{Z})^{+}$ extended to automorphisms of
$M$ by letting them act trivially on $\mathbb{U},$ the group of automorphisms
taking $(v,m,n)$ to $(v+2m\lambda,m,n-\left\langle v,\lambda\right\rangle
-m\left\langle \lambda,\lambda\right\rangle )$ for $\lambda\in L,$ and the
group of automorphisms given by reflections of norm $-2$ vectors in $R_{1}.$

\begin{lemma}
\label{Bor4}The group $\Gamma_{3}$ acts transitively on the set of vectors of
norm $-2$ in $M.$
\end{lemma}

\textbf{Proof:} The proof of Lemma \ref{Bor4} is based on the following Propositions:

\begin{proposition}
\label{Bor4a}The group $\mathcal{O}_{L}(\mathbb{Z})^{+}$ acts transitively on
the set on the set of vectors of norm $-2$ in $L.$
\end{proposition}

\textbf{Proof:} Since by definition $L=\underset{p-1}{\underbrace
{\mathbb{U}\oplus...\oplus\mathbb{U}}}\oplus\underset{q}{\underbrace
{\mathbb{E}_{8}(-1)\oplus...\oplus\mathbb{E}_{8}(-1)}}\oplus\mathbb{U}$ then
Proposition \ref{Bor4a} follows from the induction hypothesis. $\blacksquare$

\begin{proposition}
\label{Bor4b}There exists an element $\sigma\in\Gamma_{3}$ such that if
$\delta=(v,0,k)$ and $\delta^{2}=-2,$ then $\sigma(\delta)=(\mu,0,0).$
\end{proposition}

\textbf{Proof:} The condition $\left\langle \delta,\delta\right\rangle =-2$
implies that $\left\langle v,v\right\rangle =-2.$ Thus $v$ is a primitive
element in $L.$ Proposition \ref{Bor4a} implies that there exists an element
$\sigma\in\mathcal{O}_{L}(\mathbb{Z})^{+}$ such that $\sigma(\delta
)=e_{1}-e_{2}=(1,-1)\in\mathbb{U\subset}L.$ Thus easy that $\lambda=-ke_{2}\in
L$ such that $\left\langle \lambda,v\right\rangle =-k.$ Then from the
definition of the group $\Gamma_{3}$ we know that the map%
\[
\delta=(v,0,k)\rightarrow(v,0,k+\left\langle \lambda,v\right\rangle )=(v,0,0)
\]
is an automorphism. Proposition \ref{Bor4b} is proved. $\blacksquare$

\begin{proposition}
\label{Bor4c}Suppose that $\delta\in R_{0}.$ Then there exists an element
$\sigma\in\Gamma_{3}$ such that $\sigma(\delta)\in R_{1}.$
\end{proposition}

\textbf{Proof:} Proposition \ref{Bor4b} implies that without loss of
generality we may assume that $\delta\in L.$ Let $\lambda\in L,$ $\left\langle
\delta,\lambda\right\rangle \neq0$\ and $\left\langle \lambda,\lambda
\right\rangle \neq0.$ Let us consider%
\[
r=\left(  \lambda,1,-\frac{\left\langle \lambda,\lambda\right\rangle +2}%
{2}\right)  \in L\oplus\mathbb{U}=M.
\]
Clearly $\left\langle r,r\right\rangle =-2.$ Then the map $T_{r}%
(\delta)=\delta+\left\langle r,\delta\right\rangle r$ is an element of
$\Gamma_{3}$ and clearly $T_{r}(\delta)\in R_{1}.$ Proposition \ref{Bor4c} is
proved. $\blacksquare$

Combining Lemma \ref{Bor3} with Propositions \ref{Bor4a}, \ref{Bor4b} and
\ref{Bor4c} we derive Lemma \ref{Bor4}. $\blacksquare$

Lemma \ref{Bor4} implies directly Theorem \ref{Bor}. $\blacksquare$

\section{Moduli of K3 Surfaces}

\subsection{Definition of a K3 Surface}

A K3 surface is a compact, complex two dimensional manifold with the following
properties: \textbf{i. }There exists a non-zero holomorphic two form $\omega$
on X without zeroes. \textbf{ii. }$H^{1}($X,$\mathcal{O}_{\text{X}})=0.$

In \cite{Ast} and \cite{BPV}, the following topological properties are proved.
The surface X is simply connected, and the homology group $H_{2}%
(X,\mathbb{Z})$ is a torsion free abelian group of rank 22. The intersection
form $\left\langle u,v\right\rangle $ on $H_{2}(X,\mathbb{Z})$ has the
properties: \textbf{1.\ }$\left\langle u,u\right\rangle =0$ mod$(2).$
\textbf{2. }$\det\left(  \left\langle e_{i},e_{j}\right\rangle \right)  =-1$.
\textbf{3. }The symmetric form $<$ , $>$ has a signature $(3,19).$

Theorem \textbf{5} on page 54 of \cite{Se} implies that as an Euclidean
lattice $H_{2}(X,\mathbb{Z})$ is isomorphic to the K3 lattice $\Lambda_{K3},$
where $\Lambda_{K3}:=\mathbb{U}^{3}\oplus(-\mathbb{E}_{8})^{2}.$ Every K3
surface is also simply connected.

\subsection{Moduli of Marked, Algebraic and Polarized K3 surfaces}

\begin{definition}
Let $\alpha=\{\alpha_{i}\}$ be a basis of $H_{2}(X,\mathbb{Z})$ with
intersection matrix $\Lambda_{K3}.$ The pair $(X,\alpha)$ is called a marked
K3 surface. Let $l\in H^{1,1}(X,\mathbb{R})\cap H^{2}(X,\mathbb{Z})$ be the
Poincare dual class of a hyperplane section, i.e. an ample divisor. The triple
$(X,\alpha,l)$ is called a marked, polarized K3 surface. The degree of the
polarization is an integer $2d$ such that $\left\langle l,l\right\rangle
=2d>0.$
\end{definition}

\begin{definition}
The period map $\pi$ for marked K3 surfaces (X,$\alpha)$ is defined by
integrating the holomorphic two form $\omega$ along the basis $\alpha$ of
$H_{2}(X,\mathbb{Z}),$ meaning
\[
\pi(X,\alpha):=(...,\int_{\alpha_{i}}\omega,...)\in\mathbb{P}^{21}.
\]

\end{definition}

The Riemann bilinear relations hold for $\pi(X,\alpha),$ meaning
\begin{equation}
\left\langle \pi(X,\alpha),\pi(X,\alpha)\right\rangle =0\text{ }%
and\ \left\langle \pi(X,\alpha),\overline{\pi(X,\alpha)}\right\rangle >0.
\label{rbr}%
\end{equation}
Choose a primitive vector $l\in\Lambda_{K3}$ such that $\left\langle
l,l\right\rangle =2d>0.$ Let us denote $\Lambda_{K3,l}:=\left\{  v\in
\Lambda_{K3}|\left\langle l,v\right\rangle =0\right\}  .$ Then $\pi
(X,\alpha,l)\in\mathbb{P}\left(  \Lambda_{K3,l}\otimes\mathbb{C}\right)  $ and
it satisfies $\left(  \ref{rbr}\right)  .$ The set of points in $\mathbb{P}%
\left(  \Lambda_{K3,l}\otimes\mathbb{C}\right)  $ that satisfy $\left(
\ref{rbr}\right)  $ consists of two components isomorphic to the symmetric
space $\mathfrak{h}_{2,19}.$ In \cite{PS} the following Theorem was proved:

\begin{theorem}
\label{ps}The moduli space $\mathcal{M}_{K3,mpa}^{2d}$ of marked, polarized,
algebraic K3 surfaces of a fixed degree $2d$ exists and it is embedded by the
period map into $\mathfrak{h}_{2,19}$ is an open everywhere dense subset. Let
\[
\Gamma_{K3,2d}=\{\phi\in Aut^{+}(\Lambda_{K3})|\left\langle \phi
(u),\phi(u)\right\rangle =\left\langle u,u\right\rangle \text{ }and\text{
}\phi(l)=l\},
\]
where $l$ is a primitive vector such that $\left\langle l,l\right\rangle
=2d>0.$ Then the moduli space $\mathcal{M}_{K3,pa}^{2d}$ of polarized,
algebraic K3 surfaces of a fixed degree 2d is isomorphic to a Zariski open set
in the quasi-projective variety $\Gamma_{K3,2d}$%
$\backslash$%
$\mathfrak{h}_{2,19}.$
\end{theorem}

By pseudo-polarized algebraic K3 surface we understand a pair (X,$l$) where
$l$ corresponds to either ample divisor or pseudo ample divisor, which means
that for any effective divisor $D$ in X, we have $\left\langle
D,l\right\rangle \geq0.$ Mayer proved the linear system $|3l|$ defines a map:
$\phi_{|3l|}:X\rightarrow X_{1}\subset\mathbb{P}^{m}$ such that: \textbf{i.
}$X_{1}$ has singularities only double rational points. \textbf{ii. }%
$\phi_{|3l|}$ is a holomorphic birational map. Let us denote by \ $\mathcal{M}%
_{K3,ppa}^{2d}$ the moduli space of pseudo-polarized algebraic K3 surfaces of
degree 2d. From the results proved in \cite{D}, \cite{Ku}, \cite{To89} and
\cite{PS} the following Theorem follows:

\begin{theorem}
\label{pskt}The moduli space of $\mathcal{M}_{K3,ppa}^{2d}$ is isomorphic to
the locally symmetric space $\Gamma_{K3,2d}\backslash\mathfrak{h}_{2,19}.$
\end{theorem}

\subsection{Discriminant of Pseudo-Polarized K3 Surfaces}

The complement of $\mathcal{M}_{K3,mpa}^{2d}$ in $\mathfrak{h}_{2,19}$ can be
described as follow. Given a polarization class $e\in\Lambda_{K3},$ set
$T_{e}$ to be the orthogonal complement to $e$ in $\Lambda_{K3},$ i.e. $T_{e}$
is the transcendental lattice. Then we have the realization of $\mathfrak{h}%
_{2,19}$ as one of the components of $\mathfrak{h}_{2,19}\approxeq
\{u\in\mathbb{P}(T_{e}\otimes\mathbb{C})|\left\langle u,u\right\rangle =0$
$and$ $\left\langle u,\overline{u}\right\rangle >0\}.$ For each $\delta
\in\Delta(e),$ define the hyperplane $H(\delta)=\{u\in\mathbb{P}(T_{e}%
\otimes\mathbb{C})|\left\langle u,\delta\right\rangle =0\}.$ Let
\[
\mathcal{H}_{K3,2d}=\underset{\delta\in\Delta(e)}{\cup}(H(\delta
)\cap\mathfrak{h}_{2,19}).
\]
Let us define the discriminant $\mathcal{D}_{K3}^{2d}:=\Gamma_{K3,2d}%
\backslash\mathcal{H}_{K3,2d}.$ Results from \cite{Ma}, \cite{PS}, \cite{To80}
and \cite{Ku} imply that $\mathcal{D}_{K3}^{2d}$ is the complement of the
moduli space of algebraic polarized K3 surfaces $\mathcal{M}_{K3,pa}^{2d}$ in
the locally symmetric space $\Gamma_{K3,2d}\backslash\mathfrak{h}_{K3,2d},$
i.e. $\mathcal{D}_{K3}^{2d}=(\Gamma_{K3,2d}\backslash\mathfrak{h}%
_{K3,2d})\ -\mathcal{M}_{K3,pa}^{2d}.$

\subsection{Moduli of Elliptic K3 Surfaces with a Section}

\begin{definition}
\label{ell0}We will define an elliptic K3 surface X as a K3 surface such that
there exists a holomorphic map $\pi:$X$\rightarrow\mathbb{CP}^{1}$ such that
it has a section $s$.
\end{definition}

The following Theorem is a well known fact:

\begin{theorem}
\label{ell}The moduli space of elliptic pseudo polarized K3 surfaces with a
polarization class $e$ is isomorphic to $\Gamma_{ell,e}\backslash
\mathfrak{h}_{2,18},$ where $\Gamma_{ell,e}$ is defined as follows
\end{theorem}

\textbf{Proof: }Theorem \ref{ell} follows from Theorem \textbf{3.1 }proved in
\cite{dol}. $\blacksquare$

\subsection{Moduli of Einstein Metrics on K3 Surfaces}

Let $X$ be a K3 surface with a fixed $C^{\infty}$ structure. Let us consider
the set $\mathcal{M}_{\text{E}}$ of all metrics $g$ on X for which $Riccig=0$
with a volume one. We will define the moduli space $\mathfrak{M}_{\text{E}}$
of Einstein metrics as follows: $\mathfrak{M}_{\text{E}}:=\mathcal{M}%
_{\text{E}}/Diff^{+}(X)$ where $Diff^{+}(X)$ is the group of diffeomorphisms
of $X$ preserving the orientation. In \cite{KT} the following Theorem was proved:

\begin{theorem}
\label{kt}We have the following isomorphism%
\[
\mathfrak{M}_{\text{E}}\approxeq\mathcal{O}^{+}(\Lambda_{K3})\backslash
\mathbb{SO}_{0}(3,19)/\mathbb{SO}(2)\times\mathbb{SO}(19),
\]
where $\mathcal{O}^{+}(\Lambda_{K3})$ is the group of isomorphisms of the K3
lattice $\Lambda_{K3}$ which preserve the spinor norm.
\end{theorem}

\begin{lemma}
\label{a14}The moduli space of all marked elliptic K3 surfaces is an
everywhere dense subset in $\mathfrak{h}_{3,19}.$
\end{lemma}

\textbf{Proof:} We know that the set of all three dimensional vector subspaces
$E$ in $\Lambda_{K3}\otimes\mathbb{R}$ such that $E$ contains a non zero
vector in $\Lambda_{K3}$ form an every where dense subset in $\mathfrak{h}%
_{3,19}.$ Let us denote for any fixed $v\in\Lambda_{K3}$ such that
\[
\mathfrak{h}_{2,19}(v):=\left\{  E_{\tau}|\text{all oriented three dim
}E_{\tau}\subset\Lambda_{K3}\otimes\mathbb{R},\text{ }v\in E_{\tau}\text{ and
}\left\langle \text{ },\right\rangle |_{E_{\tau}}>0\right\}  .
\]
The definition of $\mathfrak{h}_{2,19}(v)$ implies that it is an open set in
the $Grass(2,21)$ and thus $\mathfrak{h}_{2,19}(v)=\mathbb{SO}_{0}%
(2,19)/\mathbb{SO}(2)\times\mathbb{SO}(19).$ It is easy to see that we can
identify $\mathfrak{h}_{2,19}(v)$ with the moduli space of marked polarized K3
surfaces with a polarization class $v.$ It is easy to see that all
$\mathfrak{h}_{2,19}(v)$ for all primitive $v\in\Lambda_{K3}$ such that
$\left\langle v,v\right\rangle >0$ form an everywhere dense subset in
$\mathfrak{h}_{3,19}.$ On the other hand it is an easy exercise to see that
all elliptic K3 surfaces with polarization vector $v$ and which have a section
form an everywhere dense subset in $\mathfrak{h}_{2,19}(v).$ For the proof of
this fact see \cite{Kod}. From here lemma \ref{a14} follows directly.
$\blacksquare$

\subsection{Moduli of K3 Surfaces with B Fields}

\begin{definition}
\label{bfield}Let $X$ be a K3 surface. A complex closed form $\omega_{X}(1,1)$
of type $(1,1)$ such that
\[%
{\displaystyle\int\limits_{X}}
\operatorname{Im}\omega(1,1)\wedge\operatorname{Im}\omega(1,1)>0
\]
will be called a B field on $X.$ The triple $(X,\alpha,\omega_{X}(1,1)),$
where $\alpha$ is a marking of the K3 surface and $\omega_{X}(1,1)$ is a
$B$-field will be called a marked K3 surface with a B-field.
\end{definition}

The moduli space of marked K3 surfaces with a B field are described by the
following Theorem:

\begin{theorem}
\label{mb}The moduli space of marked K3 surfaces with B field is isomorphic to
$\mathfrak{h}_{4,20}:=\mathbb{SO}_{0}(4,20)/\mathbb{SO}(4)\times
\mathbb{SO}(20).$See \cite{AM}.
\end{theorem}

\textbf{Proof: }The proof is based on assigning to each marked K3 surfaces $X$
with a fixed B field a four dimensional oriented subspace in $\Lambda
_{ext,K3}\otimes\mathbb{R},$ where
\[
\Lambda_{ext,K3}:=\Lambda_{K3}\oplus\mathbb{U}=\mathbb{U}\oplus\mathbb{U}%
\oplus\mathbb{U}\oplus\mathbb{E}_{8}(-1)\oplus\mathbb{E}_{8}(-1)\approxeq
H^{\ast}(\text{M,}\mathbb{Z}).
\]
Let $\omega_{X}$ be the holomorphic to form on $X.$ Then to the marked K3
surface $\left(  X;\gamma_{1},...,\gamma_{22}\right)  $ we assigne the two
dimensional oriented subspace in $\Lambda_{K3}\otimes\mathbb{R}\subset
\Lambda_{ext,K3}\otimes\mathbb{R}$ spanned by $\operatorname{Re}\omega_{X}$
and $\operatorname{Im}\omega_{X}.$ Let us assign to the B field $\omega
_{X}(1,1)$ the vector
\[
V_{X,\omega_{X}(1,1)}:=\left(  \omega_{X}(1,1),1,-\frac{\left\langle
\omega_{X}(1,1),\omega_{X}(1,1)\right\rangle }{2}\right)  \in\left(
\Lambda_{K3}\oplus\mathbb{U}\right)  \otimes\mathbb{R}.
\]
It is easy to check that
\begin{equation}
\left\langle V_{X,\omega_{X}(1,1)},V_{X,\omega_{X}(1,1)}\right\rangle =0\text{
and }\left\langle V_{X,\omega_{X}(1,1)},\overline{V_{X,\omega_{X}(1,1)}%
}\right\rangle >0. \label{RBR}%
\end{equation}
From $\left(  \ref{RBR}\right)  $ we derive that
\[
\left\langle \operatorname{Re}V_{X,\omega_{X}(1,1)},\operatorname{Re}%
V_{X,\omega_{X}(1,1)}\right\rangle =\left\langle \operatorname{Im}%
V_{X,\omega_{X}(1,1)},\operatorname{Im}V_{X,\omega_{X}(1,1)}\right\rangle >0
\]
and
\[
\left\langle \operatorname{Re}V_{X,\omega_{X}(1,1)},\operatorname{Im}%
V_{X,\omega_{X}(1,1)}\right\rangle =0.
\]
Thus the restriction of the bilinear form of $\Lambda_{ext,K3}$ to the four
dimensional subspace $E_{X,\omega_{X}(1,1)}$ spanned by $\operatorname{Re}%
\omega_{X},$ $\operatorname{Im}\omega_{X},$ $\operatorname{Re}V_{X,\omega
_{X}(1,1)},\operatorname{Im}V_{X,\omega_{X}(1,1)}$ is strictly positive. So we
constructed a map from any marked K3 surface with a B field to $\mathfrak{h}%
_{4,20}:=$ $\mathbb{SO}_{0}(4,20)/\mathbb{SO}(4)\times\mathbb{SO}(20).$
Torelli Theorem for marked K3 surfaces implies that this map is injective.

Let $E_{\tau}$ be a four dimensional positive oriented four dimensional plane
in $\left(  \Lambda_{K3}\oplus\mathbb{U}\right)  \otimes\mathbb{R}.$ Let us
consider the intersection $E_{\tau}\cap\left(  \Lambda_{K3}\otimes
\mathbb{R}\right)  .$ The orientation of $E_{\tau}$ induces an orientation on
the two dimensional positive plane $E_{\tau}\cap\left(  \Lambda_{K3}%
\otimes\mathbb{R}\right)  $ in $\Lambda_{K3}\otimes\mathbb{R}.$ From the
epimorphism of the period map we can conclude that there exists a marked K3
surface $(X,\alpha)$ such that the two dimensional space $E_{\tau}\cap\left(
\Lambda_{K3}\otimes\mathbb{R}\right)  $ is generated by $\operatorname{Re}%
\omega_{X}$ and $\operatorname{Im}\omega_{X}.$ Let us denote by $E_{X}$ the
two dimensional oriented subspace $E_{\tau}\cap\left(  \Lambda_{K3}%
\otimes\mathbb{R}\right)  $. Let $E_{X}^{\bot}$ be the orthogonal oriented
complement to $E_{X}$ in $E_{\tau}.$ Let $\omega_{1}$ and $\omega_{2}$ be two
orthonormal vectors in $E_{X}^{\bot}.$ Let us form the class of cohomology
$\Omega=\omega_{1}+\sqrt{-1}\omega_{2}\in\Lambda_{K3}\oplus\mathbb{U}.$ It is
easy to check that
\begin{equation}
\left\langle \Omega,\Omega\right\rangle =0\text{ and }\left\langle
\Omega,\overline{\Omega}\right\rangle >0. \label{RBR0}%
\end{equation}
We need the following obvious Lemma:

\begin{lemma}
\label{14}Let $e_{i}$ be the standard basis of the hyperbolic lattice
$\mathbb{U}$, i.e. $\left\langle e_{i},e_{i}\right\rangle =0$ and
$\left\langle e_{1},e_{2}\right\rangle =1.$ Let us consider $\Lambda
_{K3}\oplus\mathbb{U}.$ Then we have $\left\langle \Omega,e_{i}\right\rangle
\neq0$ for $i=1,2..$
\end{lemma}

\textbf{Proof: }Suppose that $\left\langle \Omega,e_{1}\right\rangle =0.$ Let
$\Omega_{1}:=\Omega\cap\left(  \Lambda_{K3}\otimes\mathbb{R}\right)  .$ Then
combining the assumption $\left\langle \Omega,e_{1}\right\rangle =0$ with
$\left(  \ref{RBR0}\right)  $ we get that
\begin{equation}
\left\langle \Omega_{1},\Omega_{1}\right\rangle =0\text{ and }\left\langle
\Omega_{1},\overline{\Omega_{1}}\right\rangle >0. \label{RBRa}%
\end{equation}
Thus $\left(  \ref{RBRa}\right)  $ implies that $\operatorname{Re}\Omega_{1}$
and $\operatorname{Re}\omega_{1}$ span a two dimensional positive subspace in
$\Lambda_{K3}\otimes\mathbb{R}$ which is orthogonal to the two dimensional
positive subspace $E_{X}$ in $\Lambda_{K3}\otimes\mathbb{R}.$ So their direct
sum will give a four dimensional positive subspace in $\Lambda_{K3}%
\otimes\mathbb{R}.$ This is impossible since the signature of the quadratic
form on $\Lambda_{K3}\otimes\mathbb{R}$ is $(3,19).$ Lemma \ref{14} is proved.
$\blacksquare$

Let us normalize $\Omega$ such that $\left\langle \Omega,e_{2}\right\rangle
=1.$ Then $\Omega_{1}:=\Omega\cap\left(  \Lambda_{K3}\otimes\mathbb{R}\right)
$ will satisfy the conditions $\left(  \ref{RBRa}\right)  $ which imply that
$\left\langle \operatorname{Im}\Omega_{1},\operatorname{Im}\Omega
_{1}\right\rangle >0.$ The definition of $\Omega$ implies that $\left\langle
\Omega_{1},\operatorname{Re}\Omega_{X}\right\rangle =\left\langle \Omega
_{1},\operatorname{Im}\Omega_{X}\right\rangle =0.$ Thus $\Omega_{1}$ will be a
form of type $(1,1)$ on $X$. Thus $\left\langle \operatorname{Im}\Omega
_{1},\operatorname{Im}\Omega_{1}\right\rangle >0$ implies that $\Omega_{1}$ is
a $B-$field on X. Theorem \ref{mb} is proved. $\blacksquare$

\section{Automorphic Forms on $\Gamma\backslash\mathfrak{h}_{p,q}$ and Theta
Lifts}

\subsection{Automorphic Forms of Weight -2 on $\Gamma\backslash\mathfrak{h}%
_{p,q}$}

In this paper the group $\Gamma$ will be the group of automorphisms of
$\Lambda_{K3}$ which preserve the spinor norm, i.e. $\Gamma=\mathcal{O}%
_{\Lambda_{K3}}^{+}(\mathbb{Z})$ is a subgroup of index $2$ in the group of
automorphism $\mathcal{O}_{\Lambda_{K3}}(\mathbb{Z})$ of the lattice
$\Lambda_{K3}.$ It was Donaldson who proved in \cite{D} that the mapping class
group of a K3 surface is isomorphic to $\Gamma.$ We will define the one
cocycle $\mu(\gamma,\tau)$ of the group $\mathcal{O}_{\Lambda_{K3}}%
^{+}(\mathbb{Z})$ with coefficients the non singular $3\times3$ matrices with
coefficients functions on $\mathfrak{h}_{3,19}.$

Let an element $\gamma\in\Gamma$ be represented by a matrix $(\gamma_{k,l})$
of size $(22\times22).$ Let $\tau\in\mathfrak{h}_{3,19}$ then the point $\tau$
is represented by the vectors $g_{1}(\tau),$ $g_{2}(\tau)$ and $g_{3}(\tau)$
in the fix basis $e\,_{1},...,e_{22}.$ They span a three dimensional oriented
subspace in $\mathbb{R}^{3,19}=\Lambda_{K3}\otimes\mathbb{R}$ on which the
intersection form is strictly positive. We know that the point $\tau$ can be
represented by the $3\times22$ matrix $(E_{3},\tau_{ij}),$ where $E_{3}$ is
the identity $3\times3$ matrix$.$The action of $\Gamma$ on $\mathfrak{h}%
_{3,19}$ is described as follow: Take the product of the matrices $\left(
E_{3},\tau_{i,j})\times(\gamma_{k,l})\right)  .$ It can be represented as
follows:
\begin{equation}
\left(  E_{3},\tau_{i,j})\times(\gamma_{k,l})\right)  =(\mu(\gamma
,\tau),\sigma_{\gamma,ij}(\tau)), \label{d2}%
\end{equation}
where $\mu(\gamma,\tau)$ is $3$ by $3$ matrix $A_{\tau,0}$ defined by $\left(
\ref{d2}\right)  $ and $\sigma_{\gamma,ij}(\tau)$ is some $3$ $\times19$
matrix. Theorem \ref{G4} implies that the $3\times3$ matrix $\mu(\gamma,\tau)$
has rank $3,$ i.e. $\det(\mu(\gamma,\tau))\neq0.$

\begin{definition}
\label{G6}Let $\Phi(\tau)$ be a function on $\mathfrak{h}_{3,19}$ such that it
satisfies the following functional equation:
\[
\Phi(\tau\gamma)=(\det\mu(\gamma,\tau))^{k}\Phi(\tau).
\]
Then we will call $\Phi(\tau)$ an automorphic form of weight $k.$
\end{definition}

\begin{definition}
\label{G8} Let us recall that according to Theorem \ref{G4} to each point
$\tau=(\tau_{j}^{i})\in\mathfrak{h}_{3,19},$ $1\leq j\leq3$ and $1\leq
i\leq19$ we assigned the three rows vector $g_{i}$ of the matrix ($\tau
_{j}^{i}).$ We will define the function $g(\tau)$ on $\mathfrak{h}_{3,19}$ as
follows
\begin{equation}
g(\tau):=\det(\left\langle g_{i},g_{j}\right\rangle \dot{)}. \label{f0}%
\end{equation}

\end{definition}

\begin{theorem}
\label{G9}The function $g(\tau)$ defined in $\left(  \ref{f0}\right)  $ is an
automorphic form of weight $-2.$
\end{theorem}

\textbf{Proof}: We need to compute
\[
g((\gamma(\tau))=\det(\left\langle \mu(\gamma,\tau)\times g_{i}(\gamma
(\tau)),\mu(\gamma,\tau)\times g_{j}(\gamma(\tau))\right\rangle )=?
\]
Theorem \ref{G4} and the expression of the matrix $\mu(\gamma,\tau)$ given by
$\left(  \ref{d2}\right)  $ imply
\[
g((\gamma(\tau))=\det(\mu(\gamma,\tau)\times\left\langle g_{i}(\gamma
(\tau)),g_{j}(\gamma(\tau))\right\rangle \times(\mu(\gamma,\tau))^{t})=
\]%
\[
\det(\mu(\gamma,\tau))^{2}\det(\left\langle g_{i}(\tau),g_{j}(\tau
)\right\rangle \dot{)}=\left(  \det(\mu(\gamma,\tau))^{2}\right)  \times
g(\tau).
\]
Thus we $g(\gamma(\tau))=\det(\mu(\gamma,\tau))^{2}g(\tau).$ So Theorem
\ref{G9} is proved. $\blacksquare$

\section{Regularized Determinants}

\subsection{Construction of an Automorphic Form with a Zero Set Supported by
the Discriminant Locus on $\Gamma\backslash\mathfrak{h}_{3,19}$}

The following result follows directly from the results proved in \cite{B97}.

\begin{theorem}
\label{Bor5}Let $\Lambda_{p,q}$ be an even unimodular lattice of signature
$(p,q).$ Then there exists a non zero automorphic form $\exp\left(
\Phi_{\Lambda_{p,q}}(\tau)\right)  $ such that the zero set of $\exp\left(
\Phi_{\Lambda_{p,q}}(\tau)\right)  $ coincide with the discriminant
$\mathcal{D}_{\Lambda_{p,q}}\mathcal{\subset O}_{\Lambda_{p,q}}^{+}%
(\mathbb{Z})\backslash\mathfrak{h}_{p,q}$. Moreover let $\Lambda_{p_{1},q_{1}%
}$ is an even unimodular sublattice in $\Lambda_{p,q}.$ Then
\begin{equation}
\exp\left(  \Phi_{\Lambda_{p,q}}(\tau)\right)  |_{\mathcal{O}_{\Lambda
_{p_{1},q_{1}}}^{+}(\mathbb{Z})\backslash\mathfrak{h}_{p_{1},q_{1}}}%
=\exp\left(  \Phi_{\Lambda_{p_{1},q_{1}}}(\tau)\right)  . \label{res0}%
\end{equation}

\end{theorem}

We will consider the case of K3 surfaces. We know that $\Lambda_{K3}%
=\Lambda_{3,19}.$ We will study the relations between the non zero automorphic
form $\exp\left(  \Phi_{\Lambda_{K3}}(\tau)\right)  $ and the regularized determinants.

\begin{theorem}
\label{a11}$\Delta_{B}\Phi_{\Lambda_{K3}}(\tau,\sigma)=0..$
\end{theorem}

\textbf{Proof:} Any choice of an embedding of the hyperbolic lattice
$\mathbb{U\subset}\Lambda_{K3}$ defines a totally geodesic subspace
$\mathfrak{h}_{2,18}$ into $\mathfrak{h}_{3,19}.$ This follows from Theorem
\ref{Dec}. According to the construction of the automorphic form $\exp\left(
\Phi_{\Lambda_{2,18}}(\tau)\right)  $ given in \cite{B97} it follows that
$\Phi_{\Lambda_{ell}}$ is a holomorphic function on $\mathfrak{h}_{2,18}.$
Thus we have $\Delta_{B}\Phi_{\Lambda_{ell}}=0.$ All the embeddings
$\mathfrak{h}_{2,18}\subset\mathfrak{h}_{3,19}$ corresponding to primitive
embeddings $\mathbb{U\subset}\Lambda_{K3}$ form an everywhere dense subset in
$\mathfrak{h}_{3,19}.$Since $\mathfrak{h}_{2,18}$ is a totally geodesic
subspace in $\mathfrak{h}_{3,19}$ we get that $\Delta_{B}\Phi_{\Lambda_{K3}%
}|_{\mathfrak{h}_{2,18}}=\Delta_{B}\Phi_{\Lambda_{ell}}.$ Thus the restriction
of the Bergman Laplacian applied to on $\Phi_{\Lambda_{ell}}$ is zero on an
everywhere dense subset in $\mathfrak{h}_{3,19}.$ Thus the continuous function
$\Delta_{B}\Phi_{\Lambda_{K3}}$ is zero on everywhere dense subset in
$\mathfrak{h}_{3,19}.$ From here we deduce that $\Delta_{B}\Phi_{\Lambda_{K3}%
}=0.$ Theorem \ref{a11} is proved. $\blacksquare$

\subsection{Variational Formula}

\begin{theorem}
\label{har}The function $\log\det\Delta_{KE}-\log\det\left(  \left\langle
g_{i}(\tau),g_{j}(\tau)\right\rangle \right)  $ is a harmonic function on the
moduli space $\mathfrak{M}_{\text{E}}$ of Einstein metrics of the K3 surface
with respect to the Laplacian corresponding to the Bergman metric.
\end{theorem}

\textbf{Proof: }The proof of Theorem \ref{har} is based on the following Lemmas:

\begin{lemma}
\label{har1}Let $\tau_{0}\in\mathfrak{h}_{3,19}.$ Then there exists a totally
geodesic subspace $\mathfrak{h}_{2,19}$ passing through $\tau_{0}%
\in\mathfrak{h}_{3,19}$ and its points correspond to polarized marked K3 surfaces.
\end{lemma}

\textbf{Proof: }We know that each point $\tau\in\mathfrak{h}_{3,19}$
corresponds to a three dimensional subspace $E_{\tau}\subset H^{2}%
(X,\mathbb{R})$ on which the cup product is strictly positive. Let $L\in
H^{2}(X,\mathbb{R})$ be fixed and $\left\langle L,L\right\rangle >0.$ Let us
consider the following set:%
\[
\mathfrak{h}_{L}:=\left\{  \text{three dimensional oriented positive subspaces
in }H^{2}(X,\mathbb{R})\text{ containing }L.\right\}
\]
It is easy to see that there is one to one correspondence between the two
dimensional oriented positive subspaces in the orthogonal complement $L^{\bot
}=\mathbb{R}^{2,19}$ and $\mathfrak{h}_{L}$. Thus we get that $\mathfrak{h}%
_{L}=\mathfrak{h}_{2,19}=\mathbb{SO}_{0}(2,19)/\mathbb{SO}(2)\times
\mathbb{SO}(19).$ Lemma \ref{har1} is proved. $\blacksquare$

Let us choose an orthonormal basis $e_{1},$ $e_{2}$ and $e_{3}=L$ of the three
dimensional subspace $E_{\tau_{0}}\in\mathfrak{h}_{L}.$ Lemma \ref{har1} and
Corollary \ref{G43} imply that the three dimensional subspaces $E_{\tau}$ that
correspond to $\tau\in\mathfrak{h}_{2,19}\subset\mathfrak{h}_{3,19}$ are
spanned by vectors:%
\begin{equation}
f_{1}=e_{1}+%
{\displaystyle\sum\limits_{i=1}^{19}}
\tau_{1}^{i}e_{i},\text{ }f_{2}=e_{2}+%
{\displaystyle\sum\limits_{i=1}^{19}}
{}\tau_{2}^{i}e_{i}\text{ and }f_{3}=L=e_{3}. \label{g8}%
\end{equation}

\begin{lemma}
\label{har2}In the coordinate system defined by Corollary \ref{G41} and by
$\left(  \ref{g8}\right)  $ the totally geodesic subspace is given by the
equations $\tau_{3}^{i}=0$ for $i=1,...,19.$
\end{lemma}

\textbf{Proof: }The proof follows directly from $\left(  \ref{g8}\right)  .$
$\blacksquare$

We know that $\mathfrak{h}_{2,19}$ is a complex manifold. The complex
coordinates on $\mathfrak{h}_{2,19}$ are defined as follows:%
\begin{equation}
\rho^{i}=\tau_{1}^{i}+\sqrt{-1}\tau_{2}^{i}. \label{g10}%
\end{equation}
From the epimorphism of the period map we know that $\tau_{0}$ corresponds to
a K3 surface $X_{\tau_{0}}$ and the class of cohomology of the complex two
form $e_{1}+\sqrt{-1}e_{2}$ can be identified with the class of cohomology of
the holomorphic two form $\omega_{\tau_{0}}(2,0)$ on $X_{\tau_{0}}.$ The
vector $e_{3}=L$ can be identified with the class of cohomology of the
imaginary part of a K\"{a}hler metric on $X_{\tau_{0}}.$ The subspace in
$\Lambda_{K3}\otimes\mathbb{R}$ spanned by $e_{4},...,e_{22}$ can be
identified with the primitive class of cohomology of type $(1,1),$ i.e. with
$H_{0}^{1,1}(X_{\tau_{0}},\mathbb{R})=E_{\tau_{0}}^{\perp}.$ See \cite{To80}.

We will define the Weil-Petersson metric on $\mathfrak{h}_{2,19}$ as the
restriction of the metric on $\mathfrak{h}_{3,19}$ defined by $\left(
\ref{met}\right)  .$

\begin{lemma}
\label{har3}The Weil-Petersson Metric on $\mathfrak{h}_{2,19}$ is a Hermitian metric.
\end{lemma}

\textbf{Proof:} From the expression of the Bergman metric in the coordinates
$\left(  \tau_{j}^{i}\right)  $ given by $\left(  \ref{met1}\right)  $ and
Lemmas \ref{har2} it follows that its restriction on $\mathfrak{h}_{2,19}$ is
given by
\begin{equation}
ds_{B}^{2}|_{\mathfrak{h}_{2,19}}=%
{\displaystyle\sum\limits_{i=1}^{19}}
\left(  \left(  d\tau_{1}^{i}\right)  ^{2}+\left(  d\tau_{2}^{i}\right)
^{2}\right)  +O(2). \label{b1}%
\end{equation}
Combining $\left(  \ref{g10}\right)  $ with $\left(  \ref{b1}\right)  $ we get
that%
\begin{equation}
ds_{B}^{2}|_{\mathfrak{h}_{2,19}}=%
{\displaystyle\sum\limits_{i=1}^{19}}
\left(  d\rho^{i}\right)  \otimes\left(  \overline{d\rho^{i}}\right)  +O(2).
\label{B1}%
\end{equation}
Lemma \ref{har3} is proved. $\blacksquare$

\begin{lemma}
\label{har4}Let $\tau_{0}\in\mathfrak{h}_{3,19}.$ Let $\mathfrak{h}_{2,19}$ be
the totally geodesic subspace passing through $\tau_{0}\in\mathfrak{h}_{3,19}$
and defined by the $L\in E_{\tau}$ as in Lemma \ref{har1}$.$ Then $\log$
$\det\left(  \left\langle g_{i}(\tau),g_{j}(\tau)\right\rangle \right)
|_{\mathfrak{h}_{2,19}}$ is a potential of the Weil-Petersson metric on
$\mathfrak{h}_{2,19}.$
\end{lemma}

\textbf{Proof: }Since the matrix $\left(  \left\langle g_{i}(\rho),g_{j}%
(\rho)\right\rangle \right)  |_{\mathfrak{h}_{2,19}}$ is symmetric and
\[
\left(  \left\langle g_{i}(\rho),g_{j}(\rho)\right\rangle \right)
|_{\mathfrak{h}_{2,19}}=I_{2}+\left(  h_{ij}(\rho)\right)
\]
then the following formula is true:%
\[
\log\det\left(  \left\langle g_{i}(\rho),g_{j}(\rho)\right\rangle \right)
|_{\mathfrak{h}_{2,19}}=%
{\displaystyle\sum\limits_{i=1}^{2}}
\log(1+\lambda_{i}),
\]
where $\lambda_{i}$ are the eigen values of the matrix $\left(  h_{ij}%
(\tau)\right)  .$ Thus we get%
\begin{equation}%
{\displaystyle\sum\limits_{i=1}^{2}}
\lambda_{i}=h_{11}(\tau)+h_{22}(\tau). \label{g11}%
\end{equation}
From the definition of the matrix $\left(  \left\langle g_{i}(\tau),g_{j}%
(\tau)\right\rangle \right)  |_{\mathfrak{h}_{2,19}}$ we get that
\begin{equation}
h_{11}=%
{\displaystyle\sum\limits_{i=4}^{22}}
\left(  \tau_{1}^{i}\right)  ^{2}\text{ and }h_{22}=%
{\displaystyle\sum\limits_{i=4}^{22}}
\left(  \tau_{2}^{i}\right)  ^{2}. \label{g12}%
\end{equation}
Combining $\left(  \ref{g10}\right)  ,$ $\left(  \ref{g11}\right)  $ and
$\left(  \ref{g12}\right)  $ we get that
\begin{equation}
\log\det\left(  \left\langle g_{i}(\rho),g_{j}(\rho)\right\rangle \right)
|_{\mathfrak{h}_{2,19}}=%
{\displaystyle\sum\limits_{i=4}^{22}}
|\rho^{i}|^{2}+O(3). \label{g14}%
\end{equation}
Thus we get from $\left(  \ref{g14}\right)  $ that
\begin{equation}
dd^{c}\left(  \log\det\left(  \left\langle g_{i}(\rho),g_{j}(\rho
)\right\rangle \right)  |_{\mathfrak{h}_{2,19}}\right)  =\frac{\sqrt{-1}}{2}%
{\displaystyle\sum\limits_{i=4}^{22}}
\partial\rho\wedge\overline{\partial\rho}+O(2). \label{g15}%
\end{equation}
From $\left(  \ref{g15}\right)  $ and $\left(  \ref{B1}\right)  $ we conclude
the proof of Lemma \ref{har4}. $\blacksquare$

\begin{lemma}
\label{har5}Let $\Delta_{B}$ be the Laplacian of the Bergman metric on
$\mathfrak{h}_{3,19}.$ Then the restriction of the function $\Delta_{B}\left(
\log\det\Delta_{KE}-\det\left(  \left\langle g_{i}(\tau),g_{j}(\tau
)\right\rangle \right)  \right)  $ on each totally geodesic subspace
$\mathfrak{h}_{2,19}\subset\mathfrak{h}_{3,19}$ is zero.
\end{lemma}

\textbf{Proof: }In \cite{BT05} the following Theorem was proved:

\begin{theorem}
\label{BT}Let M be a CY manifold with a polarization class $L\in H^{2}%
($M,$\mathbb{Z})\cap H^{1,1}($M,$\mathbb{R}).$ Let $\det\Delta_{(0,1)}$ be the
regularized determinant of the Laplacian corresponding to the Calabi Yau
metric corresponding to the polarization class $L$ and acting on the space of
$(0,1)$ forms. Then $dd^{c}\log\det\Delta_{(0,1)}=-\operatorname{Im}W.P..$
\end{theorem}

Combining Theorem \ref{BT} with Lemma \ref{har4} we deduce Lemma \ref{har5}.
$\blacksquare$

It is an obvious fact that the set of three dimensional positive subspaces in
$\Lambda_{K3}\otimes\mathbb{R}$ which contain a vector in $\Lambda_{K3}%
\otimes\mathbb{Q}$ form an everywhere dense subset in $\mathfrak{h}_{3,19}.$
From here it follows that we can find an everywhere dense subset of totally
geodesic subsets $\mathfrak{h}_{2,19}$ in $\mathfrak{h}_{3,19}$ on which the
continuous function
\[
\Delta_{B}\left(  \log\det\Delta_{KE}-\log\det\left(  \left\langle g_{i}%
(\tau),g_{j}(\tau)\right\rangle \right)  \right)
\]
is zero. Therefore it is zero on $\mathfrak{h}_{3,19}.$ Theorem \ref{har} is
proved. $\blacksquare$

\subsection{Relation of Regularized Determinants with Automorphic Forms}

\begin{theorem}
\label{a1}The following formula holds for the regularized determinant of the
Laplacian of the Einstein metrics $\det(\Delta_{KE})(\tau)=\det\left(
\left\langle g_{i}(\tau),g_{j}(\tau)\right\rangle \right)  \times\left\vert
\exp\left(  \Phi_{\Lambda_{K3}}(\tau)\right)  \right\vert ^{2}.$
\end{theorem}

\textbf{Proof:} According to Theorem \ref{har} the function
\[
\log\det\Delta_{KE}-\log\det\left(  \left\langle g_{i}(\tau),g_{j}%
(\tau)\right\rangle \right)
\]
is a harmonic function with respect to the Laplacian of the Bergman metric on
$\mathfrak{h}_{3,19}.$ Let us consider the function:
\[
\frac{\det\Delta_{KE}}{\det\left(  \left\langle g_{i}(\tau),g_{j}%
(\tau)\right\rangle \right)  }=\phi
\]
on $\mathfrak{h}_{3,19}.$ According to Theorem \ref{G9} the function
$\det\left(  \left\langle g_{i}(\tau),g_{j}(\tau)\right\rangle \right)  $ is
an automorphic form of weight $-2.$ Therefore the function $\phi$ is an
automorphic function of weight $2.$ In \cite{JT0} we proved that $\det
\Delta_{KE}$ is a bounded non negative function. Therefore the only zeroes of
$\det\Delta_{KE}$ can be located on the discriminant locus $\mathfrak{D}%
_{KE}.$ We know that $\left\vert \exp\left(  \Phi_{\Lambda_{K3}}(\tau)\right)
\right\vert $ is an automorphic function with a zero set on the discriminant
locus $\mathfrak{D}_{KE}.$ Since $\mathfrak{D}_{KE}$ is an irreducible in
$\mathfrak{M}_{KE},$ by taking suitable power of $\phi$ and $\left\vert
\exp\left(  \Phi_{\Lambda_{K3}}(\tau)\right)  \right\vert ,$ we may assume
that the function
\[
\frac{\left\vert \exp\left(  \Phi_{\Lambda_{K3}}(\tau)\right)  \right\vert
}{\phi}=\psi
\]
is a non zero function such $\Delta_{B}\log\psi=0.$ Thus we get a harmonic non
zero function on $\mathfrak{M}_{KE}.$

\begin{lemma}
\label{a1a} $\psi|_{\mathfrak{M}_{ell}}=const.$
\end{lemma}

\textbf{Proof: }Since $dd^{c}\left(  \log\frac{\det(\Delta_{KE})(\tau)}%
{\det\left(  \left\langle g_{i}(\tau),g_{j}(\tau)\right\rangle \right)
}\left\vert \mathfrak{M}_{ell}\right.  \right)  =0$ we can conclude that%
\[
\frac{\det(\Delta_{KE})(\tau)}{\det\left(  \left\langle g_{i}(\tau),g_{j}%
(\tau)\right\rangle \right)  }\left\vert \mathfrak{M}_{ell}\right.  =|\eta|,
\]
where $\eta$ is a holomorphic automorphic form defined up to a character
$\chi\in\Gamma_{ell}/[\Gamma_{ell},\Gamma_{ell}]$ and with a zero set
$\mathfrak{D}_{ell}.$ Since $\mathfrak{D}_{ell}$ is an irreducible divisor, we
can conclude that $\eta=\exp\left(  \Phi_{\Lambda_{ell}}(\tau)\right)  .$ Thus
since $\exp\left(  \Phi_{\Lambda_{K3}}(\tau)\right)  |_{\mathfrak{M}_{ell}%
}=\exp\left(  \Phi_{\Lambda_{ell}}(\tau)\right)  ,$ we get that $\psi
|_{\mathfrak{M}_{ell}}=conct.$ Since any two $\mathfrak{M}_{ell,1}$ and
$\mathfrak{M}_{ell,2}$ intersect. So the continuous function $\psi$ is a
constant on an everywhere dense subset in $\mathfrak{M}_{KE}$. Thus $\psi$ is
a constant. Lemma \ref{a1a} is proved. $\blacksquare$

Lemma \ref{a1a} imply Theorem \ref{a1}. $\blacksquare$

\section{Mirror Symmetry, Harvey-Moore-Borcherds Products and Counting
Problems}

\subsection{Mirror Symmetry for K3 Surfaces}

Let $(X,\alpha,\omega_{X}(1,1))$ be a marked K3 surface with a B-field. To
define the mirror of $(X,\alpha,\omega_{X}(1,1))$ we need to fix an unimodular
hyperbolic lattice $\mathbb{U}$ in $H_{2}(X,\mathbb{Z})$ with generators
$\left\{  \gamma_{0},\gamma_{1}\right\}  $ such that
\[%
{\displaystyle\int\limits_{\gamma_{0}}}
\omega_{X}\neq0\text{ }and%
{\displaystyle\int\limits_{\gamma_{1}}}
\omega_{X}\neq0.
\]
Thus we can normalize $\omega_{X}$ such that
\begin{equation}%
{\displaystyle\int\limits_{\gamma_{0}}}
\omega_{X}=1\text{ and }%
{\displaystyle\int\limits_{\gamma_{1}}}
\omega_{X}\neq0. \label{ms0}%
\end{equation}
From now on we will consider the set $(X,\alpha,\omega_{X}(1,1),\mathbb{U}),$
where $\mathbb{U}$ is a fixed sublattice in $H^{2}(X,\mathbb{Z})$ such that
the holomorphic two form satisfies $\left(  \ref{ms0}\right)  .$ Let
$\mathbb{U}^{\perp}$ be the orthogonal complement of $\mathbb{U}$ in
$H^{2}(X,\mathbb{Z}).~$Let us denote by $\mathbb{U}_{0}$ the unimodular
hyperbolic sublattice $H^{0}(X,\mathbb{Z})\oplus H^{4}(X,\mathbb{Z})$ in the
cohomology ring $H^{\ast}(X,\mathbb{Z}).$ We will assignee to the $B$-field
$\omega_{X}(1,1)$ the vector%
\[
\hat{\omega}_{X}:=\left(  \omega_{X}(1,1),1,-\frac{\omega_{X}(1,1)\wedge
\overline{\omega_{X}(1,1)}}{2}\right)
\]
in $H^{2}(X,\mathbb{Z})\oplus\mathbb{U}_{0}.$

We will need the following Theorem:

\begin{theorem}
\label{mir}There exists a marked K3 surface $Y$ with a $B$-field
$(Y,\alpha,\omega_{Y}(1,1))$ such that
\end{theorem}

\textbf{i. } \textit{If we identify }$H^{2}(Y,\mathbb{Z})$\textit{ with}
$\mathbb{U}^{\perp}\oplus\mathbb{U}_{0},$ \textit{then} $[\omega_{Y}%
]=\omega_{X}(1,1)$\textit{ in} $\left(  \mathbb{U}^{\perp}\oplus\mathbb{U}%
_{0}\right)  \otimes\mathbb{C}.$ \textbf{ii. }\textit{If we identify }%
$H^{\ast}(Y,\mathbb{Z})$ \textit{with} $H^{2}(Y,\mathbb{Z})\oplus\mathbb{U}$
\textit{then }$\omega_{Y}(1,1)=[\omega_{X}]$\textit{ in }$\left(
H^{2}(Y,\mathbb{Z})\oplus\mathbb{U}\right)  \otimes\mathbb{C},$ \textit{where
}$\omega_{X}$\textit{ is normalized as }$\left(  \ref{ms0}\right)  .$

\textbf{Proof: }Let us consider in $\left(  \mathbb{U}^{\perp}\oplus
\mathbb{U}_{0}\right)  \otimes\mathbb{C=}\Lambda_{K3}\otimes\mathbb{C}$ the
vector $\omega_{X}(1,1).$ Then direct computations show that we have
$\left\langle \hat{\omega}_{X},\hat{\omega}_{X}\right\rangle =0$ and
$\left\langle \hat{\omega}_{X},\overline{\hat{\omega}_{X}}\right\rangle >0.$
From the epimorphism of the period map for K3 surfaces proved in \cite{To80}
it follows that there exists a marked K3 surface $\,(Y,\alpha)$ with a
holomorphic two form $\omega_{Y}$ such that the class of cohomology
$[\omega_{Y}]$ is the same as the class of cohomology of $\omega_{X}(1,1).$
Next we will prove that the class of cohomology $\omega_{X}\in H^{1,1}%
(Y,\mathbb{C})$ satisfies
\[%
{\displaystyle\int\limits_{Y}}
\operatorname{Im}\omega_{X}\wedge\operatorname{Im}\omega_{X}=\left\langle
\operatorname{Im}\omega_{X},\operatorname{Im}\omega_{X}\right\rangle >0.
\]
Indeed on $X$ we have
\begin{equation}
\left\langle \omega_{X},\omega_{X}(1,1)\right\rangle =\left\langle \omega
_{X},\overline{\omega_{X}(1,1)}\right\rangle =0 \label{mir0}%
\end{equation}
since $\omega_{X}(1,1)$ is a form of type $(1,1)$ and $\omega_{X}$ is a form
of type $(2,0).$ On the other hand the form $\omega_{X}(1,1)$ with respect to
the new complex structure $Y$ on X it is a form of type $(2,0\dot{)}.$ So
$\left(  \ref{mir0}\right)  $ means that on $Y$ $\omega_{X}$ is a form of type
$(1,1\dot{)}.$On the other hand we have%
\begin{equation}%
{\displaystyle\int\limits_{X}}
\omega_{X}\wedge\overline{\omega_{X}}=2%
{\displaystyle\int\limits_{X}}
\operatorname{Im}\omega_{X}\wedge\operatorname{Im}\omega_{X}=2\left\langle
\operatorname{Im}\omega_{X},\operatorname{Im}\omega_{X}\right\rangle >0.
\label{mir1}%
\end{equation}
Thus $\left(  \ref{mir1}\right)  $ proves that $\omega_{X}$ is a B-field on
$Y.$ Theorem \ref{mir} is proved. $\blacksquare$

Now we are ready to define the mirror symmetry:

\begin{definition}
\label{MS}We will define the marked surface $(Y,\alpha,\omega_{Y}%
(1,1),\mathbb{U})$ constructed in Theorem \ref{mir} the mirror of
$(X,\alpha,\omega_{X}(1,1),\mathbb{U}_{0}).$
\end{definition}

\subsection{Mirror Symmetry and Algebraic K3 Surfaces}

Let us consider the Neron-Severi group $M=Pic(X):=H^{2}(X,\mathbb{Z})\cap
H^{1,1}(X,\mathbb{R}).$ We can characterize in another way $NS(X).$ It is dual
group in $H^{2}(X,\mathbb{Z})$ of the kernel of the functional:%
\[
\left(  \omega_{X}\right)  :H_{2}(X,\mathbb{Z})\rightarrow\mathbb{C}%
\]
defined by $\gamma\rightarrow%
{\displaystyle\int\limits_{\gamma}}
\omega_{X}.$ We define the transcendental classes of homologies $T(X)\subset$
$H^{2}(X,\mathbb{Z})$ on $X$ as follows: $T(X):=\ker\left(  \omega_{X}\right)
^{\perp}.$

We will need the following definition:

\begin{definition}
We will say that pairs $(X,M)$ $M-$marked K3 surface if $M$ is the Picard
lattice of some algebraic K3 surface together with a primitive imbedding of
$M$ into $H^{2}(X,\mathbb{Z})$.
\end{definition}

The following Theorem was proved in \cite{To93} or \cite{dol}

\begin{theorem}
\label{mmod}The moduli space $\mathfrak{M}_{M}$ of marked pairs $(X,M)$ exists
and $\mathfrak{M}_{M}\approxeq\Gamma_{M}\backslash\mathfrak{h}_{2,20-\rho},$
where $\rho=rkM$ and $\Gamma_{M}=\left\{  \phi\in Aut\Lambda_{K3}|\phi
|_{M}=id\right\}  .$
\end{theorem}

Suppose that we consider $M$ such that $\mathbb{U}$ can be embedded into
$M^{\perp}.$ According to a Theorem of Nikulin this is always possible if
$rkM=\rho\geq9.$ The construction of mirror symmetry for $M$ marked K3
surfaces $(X,\alpha,M,\mathbb{U}),$ where $\mathbb{U}\subset M^{\perp}$ was
described in \cite{To93} and \cite{dol} as follows; Let $(X,\alpha
,M,\mathbb{U})$ be an algebraic polarized K3 surface. Then Theorem \ref{mir}
implies the following Corollary:

\begin{corollary}
\label{mira}Let $(X,M,\mathbb{U},\omega_{X}(1,1))$ be $M$-marked K3 surface
such that $\mathbb{U}\subset T_{X}$ and the $B$-field $\omega_{X}(1,1)$
satisfies $\omega_{X}(1,1)|_{\mathbb{U}^{\perp}\subset T_{X}}=0.$ Then the
mirror $(Y,M_{1},\mathbb{U},\omega_{Y}(1,1))$ satisfies the following
conditions: \textbf{i. }$Pic(Y)=M_{1}=\mathbb{U}^{\perp}\subset T_{X}.$
\textbf{ii. }$T_{Y}=M\oplus\mathbb{U}\approxeq Pic(X)\oplus\mathbb{U}.$
\end{corollary}

\textbf{Proof: }Corollary \ref{mira} follows directly from Theorem \ref{mir}.
$\blacksquare$

\begin{remark}
Some interesting examples and applications of Corollary \ref{mira} were
discussed in \cite{dol}.
\end{remark}

\subsection{The Mirror Map for Marked M-K3}

Part of the mirror conjecture states that the

\begin{definition}
Let X be a K3 surface. We will define the K\"{a}hler cone of $K$(X) of X as
follows:%
\[
K\text{(X)}:=\left\{  \omega\in H^{1,1}\left(  \text{X,}\mathbb{R}\right)
|\omega=\operatorname{Im}g,\text{ and }g\text{ is a K\"{a}hler metric on
X}\right\}  .
\]

\end{definition}

We will need the characterization of the K\"{a}hler cone that is given bellow.
Denote by $\Delta(X):=\left\{  \delta\in NS(X)|\left\langle \delta
,\delta\right\rangle =-2\right\}  .$ We will need the following Lemma from
\cite{PS}:

\begin{lemma}
Let $\delta\in\Delta(X).$ Then $\delta$ or $-\delta$ can be realized as an
effective curve on X.
\end{lemma}

We will denote by $\Delta^{+}(X):=\left\{  \delta\in\Delta(X)|\delta\text{ can
be realized as an effective cure}\right\}  .$ Let us denote by $V:=\left\{
v\in H^{1,1}\left(  X,\mathbb{R}\right)  |\left\langle v,v\right\rangle
>0\right\}  .$ Since the restriction of the bilinear form on $H^{1,1}\left(
X,\mathbb{R}\right)  $ has a signature $(1,19),$ then $V$ will consists of two
components. Let us denote by $V^{+}$ the component of $V$ which contains a
K\"{a}hler class.

Each $\delta\in$ $\Delta^{+}(\Delta)$ generates a reflection $s_{\delta}$ of
$V^{+},$ where $s_{\delta}(v)=v+\left\langle v,\delta\right\rangle \delta.$
Let us denote by $\Gamma(\Delta)$ the subgroup of $\mathcal{O}_{\Lambda_{K3}%
}^{+}$ generated by $s_{\delta}.$ In \cite{To80} the following Theorem was proved:

\begin{theorem}
\label{to89}The K\"{a}hler cone $K(X)$ coincides with the fundamental domain
of the group $\Gamma(\Delta)$ in $V^{+}$ which contains a K\"{a}hler class.
\end{theorem}

\textbf{Proof: }See \cite{To80}. $\blacksquare$

\begin{remark}
\label{mirrm}According to Theorem \ref{mmod} $\mathfrak{M}_{K3,M}=\Gamma
_{M}\backslash\mathfrak{h}_{2,20-\rho}$ is the moduli space of $M$-marked K3
surfaces. Suppose that $\mathbb{U}\subset T_{X}$ is fixed and $M_{1}\subset
T_{X}$ is the orthogonal complement of $\mathbb{U}$ in $M$. Let $(Y,M_{1})$ be
some $M_{1}$ marked K3 surface defined by the primitive embedding
$M_{1}\subset T_{X}\subset\Lambda_{K3}.$ Let $\mathfrak{h}_{M_{1}}%
=M_{1}\otimes\mathbb{R}+iK(Y),$ where $K(Y)$ is the K\"{a}hler cone of $Y.$
Then according to Theorem \ref{Dec1} $\Gamma_{M}\backslash\mathfrak{h}_{M_{1}%
}\approxeq\mathfrak{M}_{K3,M}.$ Thus we have a complex analytic covering map:%
\[
\psi_{M}:\mathfrak{h}_{M_{1}}\rightarrow\Gamma_{M}\backslash\mathfrak{h}%
_{M_{1}}=\mathfrak{M}_{K3,M}.
\]
The map $\psi_{M}^{-1}$ which is multivalued is called the mirror map. It
identifies in the case described in this Remark the moduli space of $M$-marked
K3 surfaces with the complexified K\"{a}hler cone of the its mirror.
\end{remark}

\section{Applications of Mirror Symmetry}

\subsection{Counting Problems on K3}

\begin{theorem}
\label{to05}Let $X$ be an algebraic K3 surface such that Pic(X) is an
unimodular lattice. Then we have either $NS(X)=\mathbb{U\oplus E}_{8}(-1)$ or
$NS(X)=\mathbb{U\oplus E}_{8}(-1)\mathbb{\oplus E}_{8}(-1).$ Let $l\in NS(X)$
be the polarization class. Let us consider the components $V_{Enr}^{+}$ and
$V_{ell}^{+}$ of the positive cones in $\left(  \mathbb{U\oplus E}%
_{8}(-1)\right)  \otimes\mathbb{R}$ and in $\left(  \mathbb{U\oplus E}%
_{8}(-1)\mathbb{\oplus E}_{8}(-1\right)  \otimes\mathbb{R}$ which contain the
polarization vector $l.$ Let us consider the discriminant automorphic forms
$\exp\left(  \Phi_{enr}\left(  \tau\right)  \right)  $ and $\exp\left(
\Phi_{ell}\left(  \tau\right)  \right)  $ on $\left(  \mathbb{U\oplus E}%
_{8}(-1)\right)  \otimes\mathbb{R\oplus}\sqrt{-1}V_{Enr}^{+}$ and on $\left(
\left(  \mathbb{U\oplus E}_{8}(-1)\right)  \otimes\mathbb{R}\right)
\mathbb{\oplus}\sqrt{-1}V_{ell}^{+}.$ Then the restriction of the functions
$\exp\left(  \Phi_{enr}\left(  \tau\right)  \right)  $ and $\exp\left(
\Phi_{ell}\left(  \tau\right)  \right)  $ on the lines $\sqrt{-1}lt$ are
periodic. The Fourier expansions%
\[
\frac{d}{dt}\left(  \Phi_{Enr}\left(  \sqrt{-1}lt\right)  \right)  =-%
{\displaystyle\sum\limits_{n}}
a_{n}\frac{e^{-nt}}{1-e^{-nt}}%
\]
and
\begin{equation}
\text{ }\frac{d}{dt}\left(  \Phi_{ell}\left(  \sqrt{-1}lt\right)  \right)  =-%
{\displaystyle\sum\limits_{n}}
b_{n}\frac{e^{-nt}}{1-e^{-nt}}\text{ } \label{CP}%
\end{equation}
have integer coefficients $a_{n}$ and $b_{n}.$ $a_{n}$ and $b_{n\text{ }}$are
equal to the number of non singular rational curves of degree $n$ on a K3
surface $X$ with $NS(X)=\mathbb{U\oplus E}_{8}(-1)$ or $NS(X)=\mathbb{U\oplus
E}_{8}(-1)\mathbb{\oplus E}_{8}(-1).$
\end{theorem}

\textbf{Proof: }Let us fix a bases $\left\{  \gamma_{i}\right\}  $ and
$\left\{  \varepsilon_{j}\right\}  $ of $\mathbb{U\oplus E}_{8}(-1)$ and
$\mathbb{U\oplus E}_{8}(-1)\mathbb{\oplus E}_{8}(-1)$ respectively. Then we
fix the flat coordinates $\left\{  \tau^{1},...,\tau^{10}\right\}  $ and
$\left\{  \tau^{1},...,\tau^{18}\right\}  $ in the symmetric spaces
$\mathfrak{h}_{2,10}$ and $\mathfrak{h}_{2,18}$ respectfully represented as
tube domains. We will denote by $\left\langle \delta,\tau\right\rangle $ the
following expressions:%
\[
\left\langle \delta,\tau\right\rangle =%
{\displaystyle\sum\limits_{i=1}^{10}}
\left\langle \delta,\gamma_{i}\right\rangle \tau^{i}\text{ and }\left\langle
\delta,\tau\right\rangle =%
{\displaystyle\sum\limits_{i=1}^{18}}
\left\langle \delta,\varepsilon_{i}\right\rangle \tau^{i}.
\]
Then Harvey-Moore-Borcherds product formula states that there exist
automorphic forms on $\Gamma_{2,10\backslash}\mathfrak{h}_{2,10}$ or on
$\Gamma_{2,18\backslash}\mathfrak{h}_{2,18}$ which can be represented for some
large $\operatorname{Im}\tau^{i}$ as the following products.
\[
\exp\left(  \Phi_{Enr}\left(  \tau\right)  \right)  =\exp(2\pi i\left\langle
\tau,w\right\rangle
{\displaystyle\prod\limits_{\delta\in\Delta_{Enr}^{+}}}
\left(  1-\exp\left(  2\pi i%
{\displaystyle\sum\limits_{i=1}^{10}}
\left\langle \delta,\gamma_{i}\right\rangle \tau^{i}\right)  \right)
\]
and
\begin{equation}
\exp\left(  \Phi_{ell}\left(  \tau\right)  \right)  =\exp(2\pi i\left\langle
\tau,w\right\rangle
{\displaystyle\prod\limits_{\delta\in\Delta_{Enr}^{+}}}
\left(  1-\exp\left(  2\pi i%
{\displaystyle\sum\limits_{i=1}^{18}}
\left\langle \delta,\varepsilon_{i}\right\rangle \tau^{i}\right)  \right)  .
\label{cp0}%
\end{equation}
It was proved that $\exp\left(  \Phi_{Enr}\left(  \tau\right)  \right)  $ and
$\exp\left(  \Phi_{ell}\left(  \tau\right)  \right)  $ have an analytic
continuation in $\mathfrak{h}_{2,10}$ and $\mathfrak{h}_{2,18}$ and the zeroes
remain the same.

Substituting
\[%
{\displaystyle\sum\limits_{i=1}^{10}}
\gamma_{i}\tau^{i}=ilt\text{ and }%
{\displaystyle\sum\limits_{i=1}^{18}}
\varepsilon_{i}\tau^{i}=ilt\text{ }%
\]
in $\left(  \ref{cp0}\right)  $ we get
\[
\exp\left(  \Phi_{Enr}\left(  \tau\right)  \right)  =\exp(2\pi i\left\langle
\tau,w\right\rangle
{\displaystyle\prod\limits_{\delta\in\Delta_{Enr}^{+}}}
\left(  1-\exp\left(  -2\pi\left\langle \delta,l\right\rangle t\right)
\right)
\]
and
\begin{equation}
\exp\left(  \Phi_{ell}\left(  \tau\right)  \right)  =\exp(2\pi i\left\langle
\tau,w\right\rangle
{\displaystyle\prod\limits_{\delta\in\Delta_{Enr}^{+}}}
\left(  1-\exp\left(  -2\pi\left\langle \delta,l\right\rangle t\right)
\right)  . \label{cp1}%
\end{equation}
Let us split the irreducible non singular on disjoint finite sets $A_{n},$
where $A_{n}=\left\{  \delta\in\Delta^{+}|\left\langle \delta,l\right\rangle
=n\right\}  .$ Suppose that $\#A_{n}=a_{n}$ in the case of $\Lambda_{Enr\text{
}}$ and $\#A_{n}=b_{n}$ in the case $\Lambda_{ell}.$ We can rewrite $\left(
\ref{cp1}\right)  $ as follows
\[
\exp\left(  \Phi_{Enr}\left(  \tau\right)  \right)  =\exp(2\pi i\left\langle
\tau,w\right\rangle
{\displaystyle\prod\limits_{\delta\in\Delta_{Enr}^{+}}}
\left(  1-\exp\left(  -2\pi\left\langle \delta,l\right\rangle t\right)
\right)  =
\]%
\[
\exp(2\pi i\left\langle \tau,w\right\rangle
{\displaystyle\prod\limits_{n=1}}
\left(
{\displaystyle\prod\limits_{\delta\in A_{n}}}
\left(  1-\exp\left(  -2\pi nt\right)  \right)  \right)  =
\]%
\begin{equation}
\exp(2\pi i\left\langle \tau,w\right\rangle
{\displaystyle\prod\limits_{n=1}}
\left(  \left(  1-\exp\left(  -2\pi nt\right)  \right)  ^{a_{n}}\right)  .
\label{cp2}%
\end{equation}
In the same way we will get that%
\[
\exp\left(  \Phi_{ell}\left(  \tau\right)  \right)  =\exp(2\pi i\left\langle
\tau,w\right\rangle
{\displaystyle\prod\limits_{\delta\in\Delta_{Enr}^{+}}}
\left(  1-\exp\left(  -2\pi\left\langle \delta,l\right\rangle t\right)
\right)  =
\]%
\begin{equation}
\exp(2\pi i\left\langle \tau,w\right\rangle
{\displaystyle\prod\limits_{n=1}}
\left(  \left(  1-\exp\left(  -2\pi nt\right)  \right)  ^{b_{n}}\right)  .
\label{cp3}%
\end{equation}
From $\left(  \ref{cp2}\right)  $ and $\left(  \ref{cp3}\right)  $ we derive
$\left(  \ref{CP}\right)  $ and thus Theorem \ref{to05}. $\blacksquare$

\begin{remark}
We see that in the $A$ model the automorphic function $\exp\left(  \Phi
_{4,20}\left(  \tau\right)  \right)  $ restricted on the K\"{a}hler cone when
$Pic(X)$ is a unimodular lattice counts rational curves. Suppose that in the
$B$ model we represent $\mathfrak{M}_{Pic(Y)}$ as a tube domain $\mathbb{R}%
^{k}+iV^{+}$ modulo action of an arithmetic group. Suppose that the
$\operatorname{Im}\omega_{Y}\in H^{2}(Y,\mathbb{Z})\cap H^{1,1}(Y,\mathbb{R}%
).$ Then the restriction of the automorphic function $\exp\left(  \Phi
_{4,20}\left(  \tau\right)  \right)  $ on $\mathfrak{M}_{Pic(Y)}$ counts
vanishing invariant cycles $\gamma$ such that $\left\langle \gamma
,\operatorname{Im}\omega_{Y}\right\rangle =n.$
\end{remark}

\subsection{The Pluricanonical Canonical Class of the Moduli of Polarized
Algebraic K3 Surfaces}

\begin{theorem}
\label{prim}Let $l\in\Lambda_{K3}$ be a primitive vector such that
$\left\langle l,l\right\rangle =2n>0.$ Let us denote by $\left(  l\right)
^{\perp}$ be the sublattice in $\Lambda_{K3}$ orthogonal to $\mathbb{Z}l.$
Then we have%
\[
\left(  l\right)  ^{\perp}\approxeq\mathbb{Z}l^{\ast}\oplus\mathbb{U}%
^{2}\oplus\left(  -E_{8}\right)  ^{2},
\]
where $l^{\ast}$ is a primitive vector in $\Lambda_{K3}$ such that
$\left\langle l^{\ast},l^{\ast}\right\rangle =-2n<0.$
\end{theorem}

\textbf{Proof: }According to \cite{PS} the subgroup $\mathcal{O}_{\Lambda
_{K3}}^{+}$ of index two that preserve the spinor norm acts transitively on
the primitive vectors with a fixed positive self intersection. Let us fix
$\mathbb{U}$ in $\Lambda_{K3}$ with a basis $e_{0}$ and $e_{1}$ such that
$\left\langle e_{i},e_{i}\right\rangle =0$ and $\left\langle e_{1}%
,e_{2}\right\rangle =1.$ Then $l=e_{1}+ne_{2}\in\mathbb{U}$ is a primitive
vector such that $\left\langle l,l\right\rangle =2n>0.$ Let $l^{\ast}%
=e_{1}-ne_{2}\in\mathbb{U}.$ Clearly $l^{\ast}$ is a primitive vector such
that
\[
\left\langle l,l^{\ast}\right\rangle =0\text{ and }\left\langle l^{\ast
},l^{\ast}\right\rangle =-\left\langle l,l\right\rangle =-2n.
\]
Then we have
\begin{equation}
\left(  l\right)  ^{\perp}\approxeq\mathbb{Z}l^{\ast}\oplus\mathbb{U}%
\oplus\mathbb{U}\oplus\mathbb{E}_{8}(-1)\oplus\mathbb{E}_{8}(-1). \label{exp}%
\end{equation}
Theorem \ref{prim} is proved. $\blacksquare$

\begin{notation}
We will denote by $\Lambda_{K3.n}$ the lattice $\mathbb{Z}l^{\ast}%
\oplus\mathbb{U}^{2}\oplus\left(  -E_{8}\right)  ^{2}$ where $\left\langle
l^{\ast},l^{\ast}\right\rangle =-2n.$ Let $\left\{  f_{1},f_{2},g_{1}\text{
and }g_{2}\right\}  $ be a basis of $\mathbb{U}\oplus\mathbb{U}$ in $\left(
\ref{exp}\right)  $ such that $\left\langle f_{i},f_{i}\right\rangle
=\left\langle g_{i},g_{i}\right\rangle =0$ and $\left\langle f_{1}%
,f_{2}\right\rangle =\left\langle g_{1},g_{2}\right\rangle =1.$
\end{notation}

Let us consider the moduli space $\mathfrak{M}_{K3,n}$ of pseudo polarized
algebraic K3 surfaces with a polarization class $l\in\Lambda_{K3},$ where $l$
is a primitive vector in $\Lambda_{K3}$ such that $\left\langle
l,l\right\rangle =2n.$ Then according to \cite{PS} and \cite{D} we have
$\mathfrak{M}_{K3,n}=\Gamma_{n}\backslash\mathfrak{h}_{2,19},$ where
$\Gamma_{n}:=\left\{  \phi\in\mathcal{O}_{\Lambda_{K3}}^{+}|\phi(l)=l\right\}
.$ According to \cite{To80} we can define $\mathfrak{h}_{2,19}$ as one of the
open components of the quadric $\mathcal{Q}\subset\mathbb{P}(\Lambda
_{K3,n}\otimes\mathbb{C})$ defined as follows%
\[
\mathcal{Q}:=\left\{  u\in\mathbb{P}(\Lambda_{K3,n}\otimes\mathbb{C}%
)|\left\langle u,u\right\rangle =0\text{ and }\left\langle u,\overline
{u}\right\rangle >0.\right\}
\]
Let us define $\mathfrak{D}_{n}$ in $\mathfrak{M}_{K3,n}$ as follows: Let
$\lambda\in\Lambda_{K3,n},$ then%
\[
\mathcal{H}_{\lambda}:=\left\{  u\in\mathbb{P}(\Lambda_{K3,n}\otimes
\mathbb{C})|\left\langle u,\lambda\right\rangle =0.\right\}
\]
Let
\begin{equation}
\mathcal{D}_{n}=\left(
{\displaystyle\bigcup\limits_{\left\langle \delta,\delta\right\rangle
=-2\text{ \& }\delta\in\Lambda_{K3,n}}}
\left(  \mathfrak{h}_{2,19}\cap\mathcal{H}_{\delta}\right)  \right)
\cup\left(
{\displaystyle\bigcup\limits_{\phi\in\Gamma_{n}}}
\left(  \mathfrak{h}_{2,19}\cap\mathcal{H}_{\phi(l^{\ast})}\right)  \right)  .
\label{deco}%
\end{equation}
Then $\mathfrak{D}_{n}:=\Gamma_{n}\backslash\mathcal{D}_{n}.$

\begin{theorem}
\label{pc}There exists an automorphic form $\Psi_{19,n}$ on $\mathfrak{M}%
_{K3,n}=\Gamma_{n}\backslash\mathfrak{h}_{2,19}$ such that the zero set of
$\Psi_{19,n}$ is $\mathfrak{D}_{n}.$
\end{theorem}

\textbf{Proof: }According to the results of Harvey, Moore and Borcherds on we
can find an automorphic form $\left\vert \Psi_{\Lambda_{K3}}\right\vert ^{2}$
on the moduli space of Einstein metrics $\mathcal{O}_{\Lambda_{K3}}^{+}$%
$\backslash$%
$\mathfrak{h}_{3,19}$ such that its zeros are exactly on the discriminant
locus of $\mathcal{O}_{\Lambda_{K3}}^{+}$%
$\backslash$%
$\mathfrak{h}_{3,19}.$ recall that the discriminant locus on $\mathcal{O}%
_{\Lambda_{K3}}^{+}$%
$\backslash$%
$\mathfrak{h}_{3,19}$ is defined as the set of three dimensional positive
vector subspaces in $\Lambda_{K3}\otimes\mathbb{R}$ perpendicular to $\delta$
such that $\left\langle \delta,\delta\right\rangle =-2$ modulo the action of
the arithmetic group $\mathcal{O}_{\Lambda_{K3}}^{+}$. The moduli space
$\mathfrak{M}_{K3,n}=\Gamma_{n}\backslash\mathfrak{h}_{2,19}$ can be embedded
in $\mathcal{O}_{\Lambda_{K3}}^{+}$%
$\backslash$%
$\mathfrak{h}_{3,19}$ as the set of all three dimensional oriented subspaces
in $\Lambda_{K3}\otimes\mathbb{R}$ containing the polarization vector $l$
modulo the action of $\mathcal{O}_{\Lambda_{K3}}^{+}.$ The restriction of some
power of $\Psi_{\Lambda_{K3}}$ on $\mathfrak{M}_{K3,n}$ will give us an
automorphic form $\Psi_{19,n}$ on $\mathfrak{M}_{K3,n}.$ Thus we have the
following obvious fact:

\begin{remark}
\label{pc1} The zero set of $\Psi_{19,n}$ is the restriction of the zero set
of $\Psi_{\Lambda_{K3}}=\exp\left(  \Phi_{\Lambda_{3,19}}\right)  $ on
$\mathfrak{M}_{K3,n}.$
\end{remark}

Thus we need to compute the projection of the zero set of $\exp\left(
\Phi_{\Lambda_{3,19}}\right)  $ on $\Gamma^{+}\backslash\mathfrak{h}_{3,19}$
to $\mathfrak{M}_{K3,n}=\Gamma_{n}\backslash\mathfrak{h}_{2,19}.$ Theorem
\ref{pc} will follow from the following Lemma:

\begin{lemma}
\label{0} Let $\delta\in\Lambda_{K3}$ be such that $\left\langle \delta
,\delta\right\rangle =-2.$ Suppose that $\delta\notin\Lambda_{K3,n}.$ Then
there exists an automorphism $\sigma$ of the lattice $\Lambda_{K3.n}$ such
that $\sigma(\delta)\in\mathbb{U}$, i.e. $\Pr_{\mathbb{U}}\sigma
(\delta)=l^{\ast}.$
\end{lemma}

\textbf{Proof: }The proof of Lemma \ref{0} is based on the following Propositions:

\begin{proposition}
\label{2} Suppose that $\delta=me_{1}+m_{2}e_{2}+\mu_{\delta},$ $\left\langle
\delta,\delta\right\rangle =-2$ and $\delta$ satisfies $\left(  \ref{exp0}%
\right)  .$ Then there exists $\sigma\in\mathcal{O}^{+}\mathbf{(}%
\mathbb{U}^{2}\oplus\mathbb{E}_{8}(-1)^{2})$ such that
\begin{equation}
\Pr{}_{l,n}(\sigma(\delta))=n_{\delta}l^{\ast}+\mu_{\sigma(\delta)},\text{
}\mu_{\sigma(\delta)}=k_{\sigma(\delta)}\left(  f_{1}+m_{\sigma(\delta)}%
f_{2}\right)  , \label{exp2}%
\end{equation}
and
\begin{equation}
\left\Vert \mu_{\delta}\right\Vert ^{2}=\left\Vert \mu_{\sigma(\delta
)}\right\Vert ^{2} \label{exp2a}%
\end{equation}
where $k_{\sigma(\delta)}\geq1$ and $m_{\sigma(\delta)}>0.$
\end{proposition}

\textbf{Proof: }The condition that $\delta$ satisfies $\left(  \ref{exp0}%
\right)  $ implies that $\left\langle \mu_{\delta},\mu_{\delta}\right\rangle
>0.$ From the presentation of $\delta=m_{1}e_{1}+m_{2}e_{2}+\mu_{\delta}$ it
follows that have two possibilities for $\mu_{\delta}.$

\textbf{1. }$\mu_{\delta}$ is a primitive vector in $\mathbb{L}=\mathbb{U}%
\oplus\mathbb{U}\oplus\mathbb{E}_{8}(-1)\oplus\mathbb{E}_{8}(-1).$ According
to Theorem \ref{prim} all primitive vectors with a fixed positive norm form
one orbit under the action of the automorphism group. Thus there is an element
$\sigma\in\mathcal{O}_{\mathbb{L}}^{+}$ such that the primitive element
$\sigma(\mu_{\delta})$ can be presented as follows:
\[
\sigma(\mu_{\delta})=f_{1}+\frac{\left\Vert \mu_{\delta}\right\Vert ^{2}}%
{2}f_{2},
\]
where $f_{i}$ form a basis of $\mathbb{U}$ of isotropic vectors such that
$\left\langle f_{1},f_{2}\right\rangle =1.$ \textbf{2.} $\mu_{\delta}$ is not
primitive. Then the same arguments as in the first case imply $\left(
\ref{exp2}\right)  .$ Proposition \ref{2} is proved. $\blacksquare$

\begin{proposition}
\label{1} Let $\Lambda_{K3}=\mathbb{U\oplus L}$ and let $e_{1}$ and $e_{2}$ be
the isotropic generators of $\mathbb{U}.$ Let $l=e_{1}+ne_{2}\in\mathbb{U}$
and $n>0.$ Suppose that $\delta=m_{1}e_{1}+m_{2}e_{2}+\mu_{\delta}$ and
$\left\langle \delta,\delta\right\rangle =-2.$ Then there exists an element
$\sigma\in\Gamma_{n}$ such that in the representation $\sigma(\delta
)=n_{1}e_{1}+n_{2}e_{2}+\mu_{\sigma(\delta)}$
\begin{equation}
\left\langle \Pr{}_{\mathbb{U}}(\sigma(\delta)),\Pr{}_{\mathbb{U}}%
(\sigma(\delta))\right\rangle <0, \label{exp0}%
\end{equation}
where $\Pr{}_{\mathbb{U}}(\sigma(\delta)=n_{1}e_{1}+n_{2}e_{2}.$
\end{proposition}

\textbf{Proof: }Let us consider $\delta_{1}=k_{\delta_{1}}l^{\ast}+\mu
_{\delta_{1}}\in\Lambda_{K3,n},$ $\left\langle \delta_{1},\delta
_{1}\right\rangle =-2$ and $\mu_{\delta_{1}}\in\left(  l^{\ast}\right)
^{\perp}=\mathbb{L=U}\oplus\mathbb{U\oplus E}(-1)\mathbb{\oplus E}(-1).$
Clearly $\left\langle \delta_{1},\delta_{1}\right\rangle =-2$ $\ $implies
$\left\langle \mu_{\delta_{1}},\mu_{\delta_{1}}\right\rangle >0.$ Let us
consider the reflection map $\sigma(v)=r_{\delta_{1}}(v)=v+\left\langle
v,\delta_{1}\right\rangle \delta_{1},$ where $v\in\Lambda_{K3,n}.$ Then we
have $\sigma(\delta)=r_{\delta_{1}}(\delta)=\delta+\left\langle \delta
,\delta_{1}\right\rangle \delta_{1}.$ Let us compute the projection $\Pr
{}_{\mathbb{U}}(\sigma(\delta))$ of $\sigma(\delta)$ on $\mathbb{U}$ spanned
by $e_{1}$ and $e_{2}.$ Direct computations show that%
\[
\Pr{}_{\mathbb{U}}(2n\delta)=2nm_{1}e_{1}+2nm_{2}e_{2}=nm_{1}\left(
l+l^{\ast}\right)  +m_{2}\left(  l-l^{\ast}\right)  =
\]%
\[
\left(  nm_{1}+m_{2}\right)  l+\left(  nm_{1}-m_{2}\right)  l^{\ast}.
\]
Proposition \ref{2} implies that we can choose $\mu_{\delta}$ and $\mu
_{\delta_{1}}$ such that $\left\langle \mu_{\delta},\mu_{\delta_{1}%
}\right\rangle =0.$ Thus
\begin{equation}
\Pr{}_{\mathbb{U}}(\sigma(2n\delta))=\left(  nm_{1}+m_{2}\right)
l+2nk_{\delta_{1}}\left(  nm_{1}-m_{2}\right)  l^{\ast}. \label{expa}%
\end{equation}
Suppose that $nm_{1}-m_{2}\neq0.$ Then if we choose $\delta_{1}$ such that
$k_{\delta_{1}}$ is a big enouph positive number then $\left(  \ref{expa}%
\right)  $ will imply $\left(  \ref{exp0}\right)  .$

Suppose that $nm_{1}=m_{2}.$ Then $\delta=m_{1}(e_{1}+ne_{2})+\mu_{\delta}.$
Thus $\Pr{}_{\mathbb{U}}(\delta)=m_{1}l.$ Let us choose $\delta_{1}%
=k_{\delta_{1}}l^{\ast}+\mu_{\delta_{1}},$ where $\left\langle \delta
_{1},\delta_{1}\right\rangle =-2$ and $\left\langle \mu_{\delta},\mu
_{\delta_{1}}\right\rangle \neq0.$ Let us compute
\begin{equation}
r_{\delta_{1}}(\delta)=\delta+\left\langle \delta,\delta_{1}\right\rangle
\delta_{1}=m_{1}l+\left\langle \delta,\delta_{1}\right\rangle \left(
k_{\delta_{1}}l^{\ast}+\mu_{\delta_{1}}\right)  . \label{expb}%
\end{equation}
Clearly we have $\left\langle \delta,\delta_{1}\right\rangle =\left\langle
\mu_{\delta},\mu_{\delta_{1}}\right\rangle .$ Thus $\left(  \ref{expb}\right)
$ implies that
\[
r_{\delta_{1}}(\delta)=\delta_{2}=p_{1}e_{1}+p_{2}e_{2}+\left\langle
\mu_{\delta},\mu_{\delta_{1}}\right\rangle \mu_{\delta_{1}},
\]
where $np_{1}-p_{2}\neq0.$ The previous arguments imply Proposition \ref{1}.
$\blacksquare$

Let $G_{n}$ be the subgroup of $\Gamma_{n}$ generated by reflections
$r_{\kappa}(v)=v+\left\langle v,\kappa\right\rangle \kappa$ for all $\kappa
\in\Lambda_{K3,n}$ and $\left\langle \kappa,\kappa\right\rangle =-2.$ Let us
consider the orbit $G_{n}\delta$ of a fixed $\delta\in\Lambda_{K3}$ such that
$\left\langle \delta,\delta\right\rangle =-2$ and $\delta$ satisfies $\left(
\ref{exp0}\right)  .$ Let $\delta_{\min}\in\left\{  G_{n}\delta\right\}  $ be
such that
\begin{equation}
\left\langle \mu_{\delta_{\min}},\mu_{\delta_{\min}}\right\rangle =\min\left(
\underset{\delta\in\left\{  G_{n}\delta\right\}  _{+}}{\left\langle
\mu_{\delta},\mu_{\delta}\right\rangle }\right)  \geq0. \label{Min}%
\end{equation}
Proposition \ref{2} implies that without loss of generality we may suppose
that
\begin{equation}
\mu_{\delta_{\min}}=f_{1}+\frac{\left\Vert \mu_{\delta_{\min}}\right\Vert
^{2}}{2}f_{2}\text{ or }\mu_{\delta_{\min}}=k_{\delta_{\min}}\left(
f_{1}+\frac{\left\Vert \frac{\mu_{\delta_{\min}}}{k_{\delta_{\min}}%
}\right\Vert ^{2}}{2}f_{2}\right)  \label{exp3}%
\end{equation}

\begin{proposition}
\label{min}Let $\mu_{\delta_{\min}}$ be defined as $\left(  \ref{Min}\right)
.$ Then $\left\langle \mu_{\delta_{\min}},\mu_{\delta_{\min}}\right\rangle
=0.$
\end{proposition}

\textbf{Proof: }Suppose that Proposition \ref{min} is not true. Then
$\left\langle \mu_{\delta_{\min}},\mu_{\delta_{\min}}\right\rangle =\left\Vert
\mu_{\delta_{\min}}\right\Vert ^{2}>0.$ We will show that this assumption
leads to a contradiction. We can choose $\kappa\in\Lambda_{K3,n}$ such that
$\left\langle \kappa,\kappa\right\rangle =-2$ and $\kappa=k_{0}l^{\ast}%
+\mu_{\kappa},$ where $k_{0}>1.$ The relation $\left\langle \kappa
,\kappa\right\rangle =-2$ implies that $\left\Vert \mu_{\kappa}\right\Vert
^{2}=2nk_{0}^{2}-2>2n.$ Suppose that $\mu_{\delta_{\min}}$ is a primitive
element of $\mathbb{L}=\mathbb{U}\oplus\mathbb{U}\oplus\mathbb{E}%
_{8}(-1)\oplus\mathbb{E}_{8}(-1).$ Without loss of generality we can choose
\begin{equation}
\mu_{\kappa}=g_{1}+\frac{\left\Vert \mu_{\kappa}\right\Vert ^{2}}{2}%
g_{2}-f_{2}. \label{min2}%
\end{equation}
Direct computations show that $\Pr{}_{l,n}(r_{\kappa}(\delta))=k_{1}l^{\ast
}+\mu_{r_{\kappa}(\delta)},$ where
\begin{equation}
\left\Vert \mu_{r_{\kappa}(\delta)}\right\Vert ^{2}=\left\Vert \mu
_{\delta_{\min}}\right\Vert ^{2}+\left(  \left\langle \mu_{\delta_{\min}}%
,\mu_{\kappa}\right\rangle -2nk_{\delta_{\min}}k_{0}\right)  \left(
\left\Vert \mu_{\kappa}\right\Vert ^{2}+2\left\langle \mu_{\delta_{\min}}%
,\mu_{\kappa}\right\rangle \right)  \geq0. \label{Min1a}%
\end{equation}
So $\left(  \ref{exp3}\right)  $ and $\left(  \ref{min2}\right)  $ imply
\begin{equation}
\left\langle \mu_{\kappa},\mu_{\delta_{\min}}\right\rangle =\left\langle
f_{1}+\frac{\left\Vert \mu_{\delta_{\min}}\right\Vert ^{2}}{2}f_{2}%
,g_{1}+\frac{\left\Vert \mu_{\kappa}\right\Vert ^{2}}{2}g_{2}-f_{2}%
\right\rangle =-1. \label{min3}%
\end{equation}
Then from $\left(  \ref{min3}\right)  $ and $\left(  \ref{Min1a}\right)  $ we
get that $\left\Vert \mu_{\delta_{\min}}\right\Vert ^{2}>\left\Vert
\mu_{r_{\kappa}(\delta)}\right\Vert ^{2}.$ Thus we get a contradiction with
$\left\Vert \mu_{\delta_{\min}}\right\Vert ^{2}>0$ being the minimal value. So
$\left\langle \mu_{\delta_{\min}},\mu_{\delta_{\min}}\right\rangle =0.$

Suppose that $\mu_{\delta_{\min}}$ is not primitive, i.e. then $\mu
_{\delta_{\min}}=k\mu_{prim,\delta_{\min}}$ and
\[
\mu_{prim,\delta_{\min}}=f_{1}+\frac{\left\Vert \mu_{prim,\delta_{\min}%
}\right\Vert ^{2}}{2}f_{2}.
\]
Thus we get
\begin{equation}
\left\langle \mu_{\kappa},\mu_{\delta_{\min}}\right\rangle =\left\langle
k\left(  f_{1}+\frac{\left\Vert \mu_{\delta_{\min}}\right\Vert ^{2}}{2}%
f_{2}\right)  ,g_{1}+\frac{\left\Vert \mu_{\kappa}\right\Vert ^{2}}{2}%
g_{2}-f_{2}\right\rangle =-k<0. \label{min4}%
\end{equation}
Combining $\left(  \ref{Min1a}\right)  $ and $\left(  \ref{min4}\right)  $ we
get that $\left\Vert \mu_{\delta_{\min}}\right\Vert ^{2}>\left\Vert
\mu_{r_{\kappa}(\delta)}\right\Vert ^{2}.$ Thus we get a contradiction.
Proposition \ref{min} is proved. $\blacksquare$

Proposition \ref{min} implies Lemma \ref{0}. $\blacksquare$

\begin{lemma}
\label{3} The zero set of $\Psi_{19,n}$ on $\mathfrak{M}_{K3,n}$ is
$\mathfrak{D}_{n}.$
\end{lemma}

\textbf{Proof: }Let $\delta\in\Lambda_{K3}$ be such that $\left\langle
\delta,\delta\right\rangle =-2.$ Let $\Pr_{l,n}(\delta)\in\Lambda_{K3,n}$ be
the orthogonal projection of $\delta$ on $\Lambda_{K3,n}.$ If $\Pr
_{l,n}(\delta)=\delta\iff\left\langle l,\delta\right\rangle =0,$ then it
implies that the component
\[%
{\displaystyle\bigcup\limits_{\left\langle \delta,\delta\right\rangle
=-2\text{ \& }\delta\in\Lambda_{K3,n}}}
\left(  \mathfrak{h}_{2,19}\cap\mathcal{H}_{\delta}\right)
\]
in the expression $\left(  \ref{deco}\right)  $ defines one the components of
$\mathfrak{D}_{n}:=\Gamma_{n}\backslash\mathcal{D}_{n}$ corresponding to the
vectors with $-2$ norm in $\Lambda_{K3,n}.$

Suppose that $\delta\in\Lambda_{K3},$ $\left\langle \delta,\delta\right\rangle
=-2$ and $\Pr{}_{l,n}(\delta)\neq\delta.$ Lemma \ref{0} implies that we can
find $\sigma\in\Gamma_{n}$ such that $\sigma(\delta)=m_{1}e_{1}+m_{2}e_{2}.$
Thus $\Pr{}_{l,n}(\delta)=k_{\delta}l^{\ast}.$ Then
\begin{equation}
\pi\left(  \mathcal{H}_{\delta}\cap\mathfrak{h}_{2,19}\right)  =\pi\left(
H_{l^{\ast}}\cap\mathfrak{h}_{2,19}\right)  \label{min5}%
\end{equation}
where $\pi:\mathfrak{h}_{2,19}\rightarrow\Gamma_{n}\backslash\mathfrak{h}%
_{2,19}.$Thus $\left(  \ref{min5}\right)  $ implies Lemma \ref{3}.
$\blacksquare$

Theorem \ref{pc} is proved. $\blacksquare$

\end{document}